\documentclass{amsart}

\usepackage{amsmath,amsthm,amsfonts,amssymb,bm,graphicx, mathrsfs}

\usepackage
{hyperref}
\hypersetup{colorlinks=true,citecolor=blue,linkcolor=blue,urlcolor=blue,
pdfstartview=FitH }

\theoremstyle{plain}
\newtheorem{theorem}{Theorem}

\newtheorem{lemma}{Lemma}

\newtheorem{cor}[theorem]{Corollary}
\newtheorem{prop}[theorem]{Proposition}
\numberwithin{equation}{section}

\theoremstyle{definition}

\renewcommand{\geq}{\geqslant}
\newcommand{\wbar}[1]{\overline{#1}}
\renewcommand{\leq}{\leqslant}
\newcommand{\norm}[1]{\ensuremath{\left\| #1 \right\|}}
\newcommand{\e}[1]{e\left(#1\right)}
\renewcommand{\mod}{\mathrm{mod}\,}

\newcommand{\ceil}[1]{\ensuremath{\left\lceil #1 \right\rceil}}

\newcommand{\set}[1]{\ensuremath{\left\{ #1 \right\}}}

\newcommand{\changed}[1]{{\color{black} #1}}

\usepackage{calc}
\newsavebox\CBox
\newcommand\hcancel[2][0.5pt]{%
  \changed{\ifmmode\sbox\CBox{$#2$}\else\sbox\CBox{#2}\fi%
  \makebox[0pt][l]{\usebox\CBox}%
  \rule[0.5\ht\CBox-#1/2]{\wd\CBox}{#1}}}

\makeatletter
\DeclareRobustCommand\widecheck[1]{{\mathpalette\@widecheck{#1}}}
\def\@widecheck#1#2{%
    \setbox\z@\hbox{\m@th$#1#2$}%
    \setbox\tw@\hbox{\m@th$#1%
       \widehat{%
          \vrule\@width\z@\@height\ht\z@
          \vrule\@height\z@\@width\wd\z@}$}%
    \dp\tw@-\ht\z@
    \@tempdima\ht\z@ \advance\@tempdima2\ht\tw@ \divide\@tempdima\thr@@
    \setbox\tw@\hbox{%
       \raise\@tempdima\hbox{\scalebox{1}[-1]{\lower\@tempdima\box
\tw@}}}%
    {\ooalign{\box\tw@ \cr \box\z@}}}
\makeatother

\begin{document}

\author{Valentin Blomer}

\author{Jack Buttcane}
\address{Mathematisches Institut, Bunsenstr. 3-5, 37073 G\"ottingen, Germany} \email{vblomer@math.uni-goettingen.de}
\address{224 Mathematics Bldg, Buffalo, NY 14260, USA}
\email{buttcane@buffalo.edu}
 
 \title{Global decomposition of ${\rm GL}(3)$ Kloosterman sums and the spectral large sieve}

\thanks{First  author   supported  by the  Volkswagen Foundation and NSF grant 1128155 while enjoying the hospitality of the Institute for Advanced Study. The United States Government is authorized to reproduce and distribute reprints notwithstanding any copyright notation herein.}


\begin{abstract} We prove best-possible bounds for bilinear forms in Kloosterman sums for ${\rm GL}(3)$ associated with the long Weyl element. As an application we derive a best-possible spectral large sieve inequality on ${\rm GL}(3)$. 
\end{abstract}


\setcounter{tocdepth}{2}  \maketitle 

\maketitle

\section{Introduction}

\subsection{Bilinear forms in Kloosterman sums} Properties of special functions on Lie groups play an important role in the analytic theory of automorphic forms. The classical Bruggeman-Kuznetsov formula \cite{Ku} introduces integral transforms with certain Bessel kernels and their finite field analogues, classical Kloosterman sums, which are the primordial examples of algebraic exponential sums and ubiquitous in analytic number theory.  Both Bessel functions and Kloosterman sums are by now reasonably well understood. 
For groups other than ${\rm GL}(2)$, the analysis of such functions -- both over finite   rings and over the reals -- becomes very complicated. In this paper we are concerned with the special functions that arise in the long Weyl element contribution of the
${\rm GL}(3)$ Kuznetsov formula, which is typically the most interesting term in applications to spectral averages. The corresponding generalized Kloosterman sums have been worked out in full detail in \cite{BFG} and are given in explicit terms  by 
\begin{equation}\label{klooster0}
\begin{split}
&S(m_1, m_2, n_1, n_2; D_1, D_2)\\
& := \hspace{-0.4cm} \sum_{\substack{B_1, C_1 \, ({\rm mod }\, D_1)\\B_2, C_2 \, ({\rm mod }\,  D_2)\\ D_1C_2 + B_1B_2 + D_2C_1 \equiv 0 \, ({\rm mod }\, D_1D_2)\\ (B_j, C_j, D_j) = 1}} \hspace{-0.4cm} e\left(\frac{m_1B_1 + n_1(Y_1 D_2 - Z_1 B_2)}{D_1} + \frac{m_2B_2 + n_2(Y_2 D_1 - Z_2B_1)}{D_2}\right)
\end{split}
\end{equation}
for integers $n_1, n_2, m_1, m_2$ and $D_1, D_2 \in \Bbb{N}$, 
where $Y_jB_j + Z_jC_j \equiv 1 \, (\text{mod }D_j)$ for $j = 1, 2$. Whenever $(D_1, D_2) = 1$, they factorize into ordinary Kloosterman sums (\cite[Property 4.9]{BFG})
\begin{equation}\label{klooster1}
S(m_1, m_2, n_1, n_2; D_1, D_2) = S(m_1 D_2, n_1; D_1) S(n_2 D_1, m_2; D_2).
\end{equation}
On the other hand, on the diagonal $D_1 = D_2$, the Kloosterman sum \eqref{klooster0} has a very different behaviour. For instance, if $p$ is a prime not dividing $n_1n_2m_1m_2$, then
\begin{equation}\label{klooster2}
S(m_1, m_2, n_1, n_2; p, p) = p + 1
\end{equation}
is completely independent of the parameters $n_1, n_2, m_1, m_2$. Typical applications ask for a uniform treatment of these Kloosterman sums that handles simultaneously the two extreme cases \eqref{klooster1} and \eqref{klooster2} (and everything in between). One can attach an algebraic variety to the Kloosterman sum   that can be decomposed into smooth strata. This gives a corresponding local decomposition of  $S(n_1, m_2, m_1, n_2; p^r, p^s)$ with $r \leq s$ (corresponding to the $p$-adic valuation of $B_1$ and $B_2$ in \eqref{klooster0}) into pieces of the shape 
$$\sum_{x \, (\text{mod } p^t)} S(1, x; p^t) S(1, \gamma x; p^s)$$
for   certain fractional linear transformations $\gamma$ and $t \leq r/2$. This is nicely explained in \cite{St} and \cite{DF} and returns in special cases the formulae \eqref{klooster1} and \eqref{klooster2}. In this paper we give a completely explicit \emph{global} decomposition of the Kloosterman sum, at least in the case $m_1 = n_2 = 1$ which suffices for the applications we have in mind.   The formula is very complicated and will be stated in Theorem \ref{global}  in Section \ref{globaldecomp}, but roughly speaking it is an explicit interpolation between \eqref{klooster1} and \eqref{klooster2}.  As an application we are able to provide strong bounds for bilinear forms in Kloosterman sums on average over the moduli. 

\begin{theorem}\label{thm1} Let $X_1, X_2, N \geq 1$, and let $a_n$ and $b_n$ ($1 \leq n \leq N$) be two arbitrary finite sequences of complex numbers. Then
\begin{displaymath}
\begin{split}
 \mathcal{S} := & \sum_{\substack{D_1 \leq X_1\\ D_2 \leq X_2}} \Bigl| \sum_{n, m \leq N} a_n b_m S(1, m, n, 1; D_1, D_2) \Bigr| \\
&\ll (X_1X_2)^{\varepsilon} \| a\|_2 \| b \|_2 \left(X_1X_2(X_1 + X_2) + (N X_1 X_2)^{1/2} (X_1 + X_2)^{3/2} + NX_1X_2\right).
\end{split}
\end{displaymath}
\end{theorem}

As usual,  we write $\| a\|^2_2 = \sum_n |a_n|^2$.  The ``trivial''  bound, using the equivalent of Weil's bound for the Kloosterman sums \cite{St}, is $\| a \|_2 \| b \|_2 N(X_1X_2)^{3/2+\varepsilon}$.   At least if $X_1, X_2$ are not extremely unbalanced, Theorem \ref{thm1} is essentially optimal. Indeed, if $X_1 = X_2 = X$, say, then the right hand side simplifies to 
\begin{equation}\label{simple}
\| a\|_2 \| b \|_2 X^{2+\varepsilon} ( X  + N).
\end{equation}
 Considering only the cases $(D_1, D_2) = 1$ or $D_1 = D_2 = p$ and using \eqref{klooster1} and \eqref{klooster2}, one obtains a priori limitations to bounds for $\mathcal{S}$ (see \cite[Section 2]{Y2} for details) that indicate that \eqref{simple} is (up to $\varepsilon$-powers) best possible.  
 
 Our result  improves recent work of Young \cite{Y2} who obtained $\| a\|_2 \| b \|_2 X^2 (N+N^{1/2}X)^{1+\varepsilon}$ in place of \eqref{simple}. Young's work is based on the formulas \eqref{klooster1} and \eqref{klooster2}, and he uses a Fourier theoretic argument together with some rather elaborate estimates to treat the remaining cases. Our approach is more streamlined in the sense in that we use the inherent structure of the Kloosterman sum directly and avoid the detour via Fourier transform. 
 
\subsection{A spectral large sieve inequality} As an application of Theorem \ref{thm1} we provide a new and best-possible spectral large sieve inequality, the ${\rm GL}(3)$ analogue of the celebrated large sieve inequalities of Deshouillers-Iwaniec \cite{DI} that are one of the cornerstones in applications of the ${\rm GL}(2)$ Kuznetsov formula. 

 For a not necessarily cuspidal automorphic representation $\pi$ occurring in the spectral expansion of $L^2({\rm SL}_3(\Bbb{Z})\backslash \Bbb{H}_3)$ 
 let  $\lambda_{\pi}(n) $ denote the sequence of Hecke eigenvalues, normalized so that the Ramanujan conjecture predicts  $\lambda_{\pi}(n) \ll n^{\varepsilon}$,  and let $\mu_{\pi} =  \mu = (\mu_1, \mu_2, \mu_3) \in \mathfrak{a}_{\Bbb{C}}^{\ast}$ denote the spectral parameter, 
 normalized such that the Ramanujan conjecture predicts $\mu \in i \Bbb{R}^3$.  We write $$\int(\ldots) d\pi$$
 for a combined sum/integral over an orthonormal basis of spectral components of $L^2({\rm SL}_3(\Bbb{Z})\backslash \Bbb{H}_3)$, which effectively runs over Hecke-Maa{\ss} cusp forms and Eisenstein series\footnote{The constant function and the maximal Eisenstein series twisted by the constant function have no non-zero Hecke eigenvalues and therefore will not occur in \eqref{largesieve} below.} (see \cite[Theorem 10.13.1]{Go}  or \cite[Theorem 4]{Bu2} for details).  For a compact  subset $\Omega \subseteq i\mathfrak{a}^{\ast}$ we denote analogously by $\int_{\Omega} (\ldots) d\pi$ a  combined sum/integral over spectral components $\pi$ with $\mu_{\pi} \in \Omega$.    Let   
 \begin{equation}\label{weights}
\mathcal{N}(\pi) \ll \| \mu_{\pi} \|^{\varepsilon}
 \end{equation}
 be the normalizing factors ($L$-values at the edge of the critical strip) as defined in \cite[Section 3.1]{BB}. In particular, for $\pi$ cuspidal we have $\mathcal{N}(\pi)  \asymp \text{res}_{s=1} L(s, \pi \times \tilde{\pi})$. It is   conjectured that $\mathcal{N}(\pi) \gg \| \mu_{\pi}  \|^{-\varepsilon}$, but this is not known unless $\pi$ is self-dual. The best general lower bounds to date are given in \cite[Lemma 2]{Bl}. 
 
\begin{theorem}\label{thm2} Let $\Omega \subseteq i\mathfrak{a}^{\ast}$ be a compact  Weyl-group invariant subset disjoint from the Weyl chamber walls. Let $T, N \geq 1$ and $a_n$ for $N \leq n \leq 2N$ a finite sequence of complex numbers. Then
\begin{equation}\label{largesieve}
\int_{T\Omega} \frac{1}{\mathcal{N}(\pi)} \Bigl| \sum_{N \leq n \leq 2N} a_n \lambda_{\pi}(n)\Bigr|^2 d\pi \ll_{\Omega} \left(T^5  + T^2 N\right)^{1+\varepsilon} \| a \|_2^2.
\end{equation}
\end{theorem}
As we shall see below,   the bound \eqref{largesieve} is optimal (up the value of $\varepsilon$). By \eqref{weights} the weights $\mathcal{N}(\pi)^{-1}$ can be removed at no cost, if desired, but often it is convenient to include them.  Theorem \ref{thm2} improves  \cite[Theorem 3]{Bl} 
that had a factor $N^2$  in the second term. A similar argument can also be used to improve Young's local version of the large sieve \cite[Theorem 1.1]{Y2} by a factor $N^{1/2}$.  We leave the details for the   (analytically   simpler) local version to the reader. \\

At first sight it may be surprising that  Theorem \ref{thm2} is   optimal. From general principles of the large sieve (\cite[Section 7]{IK}) it is clear that one cannot do better than 
\begin{equation}\label{perfect}
 ( \text{number of harmonics } +  \text{ length of summation}) \| a \|_2^2 \asymp  (T^5 + N)   \| a \|_2^2
\end{equation}  
   on the right hand side of \eqref{largesieve}. 
    However, in the situation of Theorem \ref{thm2} we have the following somewhat unexpected lower bound.
    
 \begin{prop}\label{lowerbound} Suppose that $\Omega$ has non-empty interior and that $N \geq T^{3+\delta}$ for some $\delta > 0$ and some sufficiently large $T$. Then there exists a sequence $a_n$  such that
  $$ \int_{T\Omega} \frac{1}{\mathcal{N}(\pi)} \Bigl| \sum_{N \leq n \leq 2N} a_n \lambda_{\pi}(n)\Bigr|^2 d\pi \gg  T^{2}  N  \| a \|_2^2.$$
 \end{prop}

 Natural families, for which ``perfect'' large sieve inequalities of the type \eqref{perfect} fail, are rather rare, the only other prominent example being the family of Fourier coefficients of cusp forms for $\Gamma_1(N)$ \cite{IL}. Our example produces a new family of this kind, and it will be clear from the proof of Proposition \ref{lowerbound}  in Section \ref{proof}  that this generalizes easily to higher rank: a spectral large sieve of the type of Theorem \ref{thm2} on ${\rm GL}(r)$, $r \geq 2$, cannot be stronger than $$(T^{\dim \Bbb{H}_r} + T^{\dim \Bbb{H}_{r-1}}N) \| a \|_2^2,$$
 where $\dim \Bbb{H}_r = r(r+1)/2 - 1$ is the dimension of the generalized upper half plane.     \\ 

As an application of Theorem \ref{thm2} we present the following uniform Lindel\"of-on-average bound for a second moment of a degree 6 family of $L$-functions. 
\begin{cor}\label{cor3}  Let $\Omega$ be as in Theorem \ref{thm2} and suppose in addition that $\Omega$ is disjoint from the set $\mu_1\mu_2\mu_3 = 0$. For $T \geq 1$ let $t \in [-T^{1-\varepsilon}, T^{1-\varepsilon}]$ and let $f$ be an even Hecke-Maa{\ss} cusp form  for ${\rm SL}_2(\Bbb{Z})$ with spectral parameter $\tau \leq T^{1-\varepsilon}$. Then
$$\int_{T\Omega} |L(1/2 + it, \pi \times f)|^2 d\pi \ll  T^{5+\varepsilon}$$
where the implied constant depends only on $\varepsilon$ and is independent of $t$ and $f$. 
\end{cor}

Theorem \ref{thm2} is proved by opening the square and applying the Kuznetsov formula. All but the long Weyl element contribution are simple to estimate. For the latter we couple a variation of Theorem \ref{thm1} with a   hybrid large sieve for ${\rm GL}(1)$ harmonics. At this point also the archimedean integral transforms in the Kuznetsov formula enter the picture. We devote Section \ref{analysis} to the very delicate analysis of these special functions, where we investigate  averages of the  Whittaker transforms   over a large (smooth) region $T \Omega$ as in Theorem \ref{thm2}.   As is to be expected, the rather complicated asymptotic behaviour of these transforms stabilizes under such an average, but the result is more complex than one might expect at first sight. As for the Kloosterman sums \eqref{klooster0}, there is a clear dichotomy between values close to the diagonal and away from the diagonal. A formal stationary phase analysis shows an oscillation of the type 
$$e\left( \pm \big|y_1^{1/3} \pm y_2^{1/3}\big|^{3/2}\right).$$
The algebraic shape of the oscillating factor is rather characteristic for special functions on ${\rm GL}(3)$, cf.\ e.g.\ \cite{BuHu} for a somewhat related situation. To keep the paper at reasonable length, we prove only as much as is needed for the application at hand and perform a few short-cuts to avoid multi-dimensional stationary phase analysis.  \\

Two comments on \textbf{notation}. As is customary, we use $\varepsilon$-convention most of the time. Certain arguments require, however, a careful treatment of epsilon-powers, so that locally in some places exponents like $2\varepsilon$ or $\varepsilon/2$ will occur. 
For two quantities $A, B > 0$   we write $A \asymp B$ to mean that there are positive constants $c_1$, $c_2$ such that $c_1A  \leq B \leq c_2A$.

\section{Global decomposition of the Kloosterman sum}\label{globaldecomp}
 
 We recall the definition  \eqref{klooster0} of the Kloosterman sum associated with the long Weyl element of the group ${\rm GL}(3)$. 

\subsection{Preliminary decompositions}
In this subsection we assume
\begin{equation}\label{assume} 
 p\mid D_1 \Leftrightarrow p\mid D_2 \quad \text{and} \quad  1 <  D_1\mid D_2.
 \end{equation}
We start by considering various subsums of the Kloosterman sum \eqref{klooster0}. 

Let $$S_1 = S_1(m_1, m_2, n_1, n_2; D_1, D_2)$$ be the sum subject to the extra condition $(B_1, D_1) = 1$. The summation conditions imply $D_1 \mid B_2$, so we substitute $B_2 \mapsto B_2 D_1$, and extend both $B_j$ sums modulo $D_2$.
We may choose $Z_1=Y_2=0$, $Y_1 = \wbar{B_1}$, $Z_2 = \wbar{C_2}$, so
\begin{align*}
	S_1  = \frac{1}{\frac{D_2}{D_1} D_1} \sum_{\substack{B_1, B_2 \, (\mod{D_2})\\(B_1,D_1)=1}} \sum_{\substack{C_1\, (\mod{D_1})\\(B_1 B_2+C_1(D_2/D_1),D_2)=1}} \e{\frac{m_1 B_1}{D_1}+\frac{m_2 D_1 B_2+n_2 B_1 \wbar{\left(B_1 B_2+C_1(D_2/D_1)\right)}}{D_2}}.
\end{align*}
Now substitute $C_1 \mapsto B_1 C_1$ and then $B_2 \mapsto B_2-C_1(D_2/D_1)$, so $(B_2,D_2)=1$ and the $C_1$ sum drops out, giving
\begin{align*}
	S_1  =& \frac{D_1}{D_2} \sum_{\substack{B_1, B_2 \, (\mod{D_2})\\(B_1,D_1)=(B_2,D_2)=1}} \e{\frac{m_1 B_1}{D_1}+\frac{m_2 D_1 B_2+n_2 \wbar{B_2}}{D_2}} =  S(0, m_1; D_1) S(m_2 D_1, n_2; D_2).
\end{align*}

Similarly, defining $S_2 = S_2(m_1, m_2, n_1, n_2; D_1, D_2)$ to be the sums with $(B_2, D_2) = 1$, we obtain
$$S_2 =  S(0,n_1; D_1) S(m_2, n_2 D_1; D_2).$$

Next, let $D_1^\dagger$ and $D_2^{\dagger}$ denote the squarefree kernels of $D_1$ and $D_2$. For $$D_1^\dagger \mid E_1 \mid D_1, \quad D_1^\dagger \mid E_2 \mid D_2, \quad   D_1 \mid E_1 E_2$$
let $$S_3'(E_1, E_2) = S_3'(E_1, E_2; m_1, m_2, n_1, n_2;D_1, D_2)$$ be the subsum with the conditions $E_1 \mid B_1$, $E_2 \mid B_2$ and $(B_j/E_j,D_j)=1$ for $j= 1, 2$. We replace $B_j \mapsto B_j E_j$ and note that the conditions imply $(C_j,D_j)=1$, so we choose $Y_j=0$ and $Z_j = \wbar{C_j}$.
We now extend the new $B_j$  sums to be mod $D_2$,  
so that
\begin{align*}
	S_3'( E_1, E_2) = \frac{C_{D_1, E_1} C_{D_2, E_2}}{ (D_2 E_1/D_1)E_2} &\sum_{\substack{B_1, B_2\, (\mod{D_2})\\(B_1,D_1)=(B_2,D_2)=1}} \sum_{\substack{C_1 \, (\mod{D_1})\\(C_1,D_1)=1\\(B_1 B_2 \frac{E_1 E_2}{D_1}+C_1\frac{D_2}{D_1},D_2)=1}} \e{\frac{m_1 E_1 B_1-n_1\wbar{C_1}B_2 E_2}{D_1}}  \\
& \times   \e{\frac{m_2 B_2 E_2+n_2 E_1 B_1 \wbar{\left(B_1 B_2 (E_1 E_2/D_1)+C_1(D_2/D_1)\right)}}{D_2}}.
\end{align*}
Here 
\begin{equation}\label{factor}
C_{D_1,E_1} := \prod_{v_p(D_1)=v_p(E_1)} \left(1-\frac{1}{p}\right)^{-1}
\end{equation}
(where $v_p$ is the usual $p$-adic valuation)  occurs in case $E_1$ contains the entirety of a prime-power factor of $D_1$, in which case extending the $B_1$ sum required skipping the terms with $p\mid B_1$ (to maintain $(B_1, D_1)=1$). Moreover, we have removed $C_2$ from the summation by the constraint equation $C_2+B_1 B_2 (E_1 E_2/D_1)+C_1 (D_2/D_1) \equiv 0\, (\mod{D_2})$.

Now we may substitute $C_1 \mapsto -B_1 B_2 C_1$, and write
\begin{align*}
	C_1^* = \frac{E_1 E_2}{D_1}-C_1\frac{D_2}{D_1},
\end{align*}
giving
\begin{align*}
	S_3' ( E_1, E_2)=& \frac{C_{D_1, E_1} C_{D_2, E_2}}{E_2 (D_2 E_1/D_1)} \sum_{\substack{B_1, B_2\, (\mod{D_2})\\(B_1,D_1)=(B_2,D_2)=1}} \sum_{\substack{C_1 \, (\mod{D_1})\\(C_1,D_1)=(C_1^*,D_2)=1}} \e{\frac{m_1 E_1 B_1+n_1\wbar{C_1} E_2 B_1}{D_1}}\\
& \quad\quad\quad\quad	\times  \e{\frac{m_2 B_2 E_2+n_2 E_1 B_2 \wbar{C_1^*}}{D_2}} \\
	=& \frac{C_{D_1, E_1} C_{D_2, E_2}}{E_1 E_2} \sum_{\substack{C_1 \, (\mod{D_1})\\(C_1,D_1)=(C_1^*,D_2)=1}} S(m_1 E_1, n_1 E_2 \wbar{C_1}; D_1) S(m_2 E_2, n_2 E_1 \wbar{C_1^*};D_2).
\end{align*}
Notice that the $C_1$ sum is empty unless $(E_1 E_2/D_1, D_2/D_1) = 1$.
In particular, we must have $E_2 \mid D_1$.

Now let $S_3 = S_3(m_1, m_2, n_1, n_2; D_1, D_2)$ be the subsum with $D_1^{\dagger} \mid B_1$, $D_1^{\dagger} \mid B_2$. Then we conclude from the above that
\begin{displaymath}
\begin{split}
S_3 &=  \sum_{\substack{D_1^\dagger \mid E_j \mid D_1 \\ (E_1 E_2, D_2)=D_1}}  S_3'( E_1, E_2)\\
& = \sum_{\substack{D_1^\dagger \mid E_j \mid D_1 \\ (E_1 E_2, D_2)=D_1}} \frac{C_{D_1, E_1} C_{D_2, E_2}}{E_1 E_2} \sum_{\substack{C_1 \,(\mod{D_1})\\(C_1,D_1)=(C_1^*,D_2)=1}} S(m_1 E_1, n_1 E_2 \wbar{C_1}; D_1) S(m_2 E_2, n_2 E_1 \wbar{C_1^*};D_2).
\end{split}
\end{displaymath}

We can collect the decompositions of $S_1, S_2, S_3$ as follows.  For any $e_1 e_2 e_3 = D_1$, $f_1 f_2 f_3=D_2$ (where $D_1, D_2$ still satisfy \eqref{assume}) with  
\begin{equation}\label{prime}
  e_i \text{ pairwise coprime,} \quad  p \mid f_i \Leftrightarrow p\mid e_i,
  \end{equation}  
  we let $S_4 = S_4(m_1, m_2, n_1, n_2; D_1, D_2)$ be the sum with conditions $(B_1,e_1)=1$, $(B_2, f_2)=1$, $e_3^\dagger \mid B_1$, $e_3^{\dagger} \mid B_2$.  By the Chinese Remainder Theorem \cite[Property 4.7]{BFG}, we have
\begin{align*}
	 S_4 & = S_1(\wbar{e_2 e_3}^2 f_2 f_3 m_1, \wbar{f_2 f_3}^2 e_2 e_3 m_2, n_1, n_2; e_1, f_1) \\
	& \qquad \times S_2(\wbar{e_1 e_3}^2 f_1 f_3 m_1, \wbar{f_1 f_3}^2 e_1 e_3 m_2, n_1, n_2; e_2, f_2) \\
	& \qquad \times S_3(\wbar{e_1 e_2}^2 f_1 f_2 m_1, \wbar{f_1 f_2}^2 e_1 e_2 m_2, n_1, n_2; e_3, f_3) \\
	&= S(0, m_1; e_1) S(\wbar{f_2 f_3}^2 D_1 m_2, n_2; f_1)   S(0,n_1; e_2) S(\wbar{f_1 f_3}^2 m_2, D_1 n_2; f_2)\sum_{\substack{e_3^\dagger | E_j | e_3 \\ (E_1 E_2, f_3)=e_3}} \frac{C_{e_3, E_1} C_{f_3, E_2}}{E_1 E_2} \\
	& \qquad \times  \sum_{\substack{C_1 \, (\mod{e_3})\\(C_1,e_3)=(C_1^*,f_3)=1}} S(\wbar{e_1 e_2}^2 f_1 f_2 m_1 E_1, n_1 E_2 \wbar{C_1}; e_3) S(\wbar{f_1 f_2}^2 e_1 e_2 m_2 E_2, n_2 E_1 \wbar{C_1^*};f_3),
\end{align*}
where of course the meaning of $C_1^{\ast}$ is now
$$C_1^* = \frac{E_1 E_2}{e_3}-C_1\frac{f_3}{e_3}.$$

Now we finally obtain
$$ S(m_1,m_2,n_1,n_2; D_1, D_2) = \mathop{{\sum}'}\limits_{\substack{e_1 e_2 e_3=D_1\\f_1 f_2 f_3=D_2}} S_4(m_1,m_2,n_1,n_2; e_1 e_2 e_3, f_1 f_2 f_3), $$
 where the prime indicates the conditions \eqref{prime}. Taking $E_3=(E_1,E_2)$, we may write  the right hand side of the preceding equation as
 \begin{align*}
	&\mathop{{\sum}'}\limits_{\substack{e_1 e_2 e_3=D_1\\f_1 f_2 f_3=D_2}} \sum_{\substack{e_3^\dagger \mid E_3 \mid e_3, E_1 E_3 \mid e_3, E_2 E_3 \mid e_3 \\ (E_1 E_2 E_3^2, f_3)=e_3 \\ (E_1,E_2)=1}} \sum_{\substack{C_1 \, (\mod{e_3/E_3})\\(C_1,e_3/E_3)=(C_1^*,f_3/E_3)=1}} \frac{E_3 C_{e_3, E_1 E_3} C_{f_3, E_2 E_3}}{E_1 E_2 (C_{e_3, E_3})^2 C_{f_3, E_3}} \\
	& \qquad \times  S(0, m_1; e_1) S(0,n_1; e_2) S(\wbar{f_2 f_3}^2 D_1 m_2, n_2; f_1) S(\wbar{f_1 f_3}^2 m_2, D_1 n_2; f_2) \\
	& \qquad \times S(\wbar{e_1 e_2}^2 f_1 f_2 m_1 E_1, n_1 E_2 \wbar{C_1}; e_3/E_3) S(\wbar{f_1 f_2}^2 e_1 e_2 m_2 E_2, n_2 E_1 \wbar{C_1^*};f_3/E_3),
\end{align*}
where with the current choice of variables   $$C_1^* = \frac{E_1 E_2 E_3^2}{e_3}-C_1\frac{f_3}{e_3}. $$
We give two special cases. First suppose $v_p(D_2) > v_p(D_1)$ for each $p$, then the condition $(E_1 E_2 E_3^2, f_3)=e_3$ becomes $E_1 E_2 E_3^2=e_3$ and $v_p(E_3),v_p(E_1 E_3), v_p(E_2 E_3)<v_p(e_3)<v_p(f_3)$, so we may write
\begin{align*}
	& S(m_1,m_2,n_1,n_2; D_1, D_2) \\
	&= \mathop{{\sum}'}\limits_{\substack{e_1 e_2 e_3=D_1\\f_1 f_2 f_3=D_2\\e_1\mid n_2, e_2\mid m_2}} \sum_{\substack{E_1 E_2 E_3^2=e_3 \\ e_3^\dagger \mid  E_3, (E_1,E_2)=1}} \sum_{\substack{C_1 \, (\mod{e_3/E_3})\\(C_1,e_3/E_3)=(C_1^*,f_3/E_3)=1}} \frac{e_1 e_2 E_3}{E_1 E_2} S(0, m_1; e_1) S(0,n_1; e_2) \\
	& \qquad \times S(\wbar{f_2 f_3}^2 (D_1/e_1) m_2, (n_2/e_1); f_1/e_1) S(\wbar{f_1 f_3}^2 (m_2/e_2), (D_1/e_2) n_2; f_2/e_2) \\
	& \qquad \times S(\wbar{e_1 e_2}^2 f_1 f_2 m_1 E_1, n_1 E_2 \wbar{C_1}; e_3/E_3) S(\wbar{f_1 f_2}^2 e_1 e_2 m_2 E_2, n_2 E_1 \wbar{C_1^*};f_3/E_3).
\end{align*}

On the other hand, when $D_1=D_2 $, the decomposition becomes
\begin{align*}
	& S(m_1,m_2,n_1,n_2; D, D) \\
	&= \mathop{{\sum}'}\limits_{\substack{e_1 e_2 e_3=D}} \sum_{\substack{e_3^\dagger \mid E_3 \mid e_3, E_1 E_3 \mid e_3\\  E_2 E_3 \mid e_3,  e_3\mid E_1 E_2 E_3^2 \\ (E_1,E_2)=1}} \sum_{\substack{C_1 (\mod{e_3/E_3})\\(C_1,e_3/E_3)=(C_1^*,e_3/E_3)=1}} \frac{E_3 C_{e_3, E_1 E_3} C_{e_3, E_2 E_3}}{E_1 E_2 (C_{e_3, E_3})^3} S(0, m_1; e_1) S(0,n_1; e_2) \\
	& \qquad \times S(0, n_2; e_1) S(m_2, 0; e_2) S(\wbar{e_1 e_2} m_1 E_1, n_1 E_2 \wbar{C_1}; e_3/E_3) S(\wbar{e_1 e_2} m_2 E_2, n_2 E_1 \wbar{C_1^*};e_3/E_3).
\end{align*}

\subsection{The full decomposition}
We now remove the assumptions \eqref{assume} on $D_1$ and $D_2$. 
Let $F_j G_j H_j I = D_j$ where $F_j$, $G_j$, $H_j$, and $I$ are coprime for each $j=1,2$ and we have
\begin{displaymath}
\begin{split}
& (F_1,F_2)=1, \quad p\mid G_1\Leftrightarrow p\mid G_2, \quad p\mid H_1\Leftrightarrow p\mid H_2, \\
 & p \mid G_1 \Rightarrow v_p(G_2) > v_p(G_1),\quad  p\mid H_1 \Rightarrow v_p(H_1) > v_p(H_2).
 \end{split}
 \end{displaymath}
We now specialize to the situation where $m_1 = n_2 = 1$. We factor
\begin{align*}
	 S(1, m, n, 1; D_1, D_2) &= S((\wbar{D_1/F_1})^2 (D_2/F_2), (\wbar{D_2/F_2})^2 (D_1/F_1) m, n, 1; F_1, F_2) \\
	&\qquad S((\wbar{D_1/G_1})^2 (D_2/G_2), (\wbar{D_2/G_2})^2 (D_1/G_1) m, n, 1; G_1, G_2) \\
	&\qquad S((\wbar{D_2/H_2})^2 (D_1/H_1) m, (\wbar{D_1/H_1})^2 (D_2/H_2), 1, n; H_2, H_1) \\
	&\qquad S((\wbar{D_1/I})^2 (D_2/I), (\wbar{D_2/I})^2 (D_1/I) m, n, 1; I, I),
\end{align*}
and apply the above decompositions to obtain 
\begin{equation}\label{multiline}
\begin{split}
	& S(1, m, n, 1; D_1, D_2) = S((\wbar{D_1/F_1})^2 D_2, n; F_1) S((\wbar{D_2/F_2})^2 D_1, m; F_2) \\
	&   \times \mathop{{\sum}'}\limits_{\substack{e_2 e_3=G_1\\f_2 f_3=G_2\\e_2|m}} \sum_{\substack{E_1 E_2 E_3^2=e_3 \\ e_3^\dagger \mid E_3, (E_1,E_2)=1}} \sum_{\substack{C_1 \, (\mod{e_3/E_3})\\(C_1,e_3/E_3)=(C_1^*,f_3/E_3)=1}} \frac{e_2 E_3}{E_1 E_2} S(0,n; e_2) S\left(\Bigl(\wbar{\frac{D_2}{f_2}}\Bigr)^2 \frac{D_1}{e_2} \frac{m}{e_2}, 1; \frac{f_2}{e_2}\right) \\
	&  \qquad \times S\left(E_1 \Bigl(\wbar{\frac{D_1}{e_3}}\Bigr)^2 \frac{D_2}{f_3}, n E_2 \wbar{C_1}; \frac{e_3}{E_3}\right) S\left(\Bigl(\wbar{\frac{D_2}{f_3}}\Bigr)^2 \frac{D_1}{e_3} E_2 m, E_1 \wbar{C_1^*};\frac{f_3}{E_3}\right) \\
	&  \times \mathop{{\sum}'}\limits_{\substack{e_4 e_6=H_2\\f_4 f_6=H_1\\e_4\mid n}} \sum_{\substack{E_4 E_5 E_6^2=e_6 \\ e_6^\dagger \mid E_6, (E_4,E_5)=1}} \sum_{\substack{C_2 \, (\mod{e_6/E_6})\\(C_2,e_6/E_6)=(C_2^*,f_6/E_6)=1}} \frac{e_4 E_6}{E_4 E_5} S(0, m; e_4) S\left(\Bigl(\wbar{\frac{D_1}{f_4}}\Bigr)^2 \frac{D_2}{e_4}, \frac{n}{e_4}; \frac{f_4}{e_4}\right) \\
	& \qquad \times S\left(E_4 \Bigl(\wbar{\frac{D_2}{e_6}}\Bigr)^2 \frac{D_1}{f_6} m, E_5 \wbar{C_2}; \frac{e_6}{E_6}\right) S\left(E_5 \Bigl(\wbar{\frac{D_1}{f_6}}\Bigr)^2 \frac{D_2}{e_6}, n E_4 \wbar{C_2^*};\frac{f_6}{E_6}\right) \\
	&   \times \mathop{{\sum}'}\limits_{\substack{e_7 e_8 e_9=I}} \sum_{\substack{e_9^\dagger \mid E_9 \mid e_9, E_7 E_9 \mid e_9\\ E_8 E_9 \mid e_9, e_9\mid E_7 E_8 E_9^2 \\ (E_7,E_8)=1}} \sum_{\substack{C_3 \, (\mod{e_9/E_9})\\(C_3,e_9/E_9)=(C_3^*,e_9/E_9)=1}} \frac{E_9 C_{e_9, E_7 E_9} C_{e_9, E_8 E_9}}{E_7 E_8 (C_{e_9, E_9})^3} S(0, 1; e_7) S(0,n; e_8) \\
	& \qquad \times S(0, 1; e_7) S(m, 0; e_8) S\left(E_7 \Bigl(\wbar{\frac{D_1}{e_9}}\Bigr)^2 \frac{D_2}{e_9}, n E_8 \wbar{C_3}; \frac{e_9}{E_9}\right) S\left(E_8 \Bigl(\wbar{\frac{D_2}{e_9}}\Bigr)^2 \frac{D_1}{e_9} m, E_7 \wbar{C_3^*};\frac{e_9}{E_9}\right),
	\end{split}
\end{equation}
where $$C_1^* = 1-C_1(f_3/e_3), \quad C_2^* = 1-C_2(f_6/e_6), \quad C_3^* = (E_7 E_8 E_9^2/e_9)-C_3.$$  
 
We continue with further manipulations. We need to hold the Ramanujan sums separate as they are very small on the factors which are coprime to m and n. We also need to work out the $E_1E_2E_4E_5E_7E_8$ factors as they will artificially increase the modulus (and hence the waste in the $F_1$ sum) below. We consider a Kloosterman  sum of the form $S(a n, b\alpha; c)$ with $ab \mid c$ and $(a,b)=(\alpha,c)=1$, so that
\begin{align*}
	S(a n, b\alpha; c) =& \sum_{\substack{x\, (\mod{c/a''})\\(x,c/a'')=1}} \e{\frac{an\wbar{x}+b\alpha x}{c}} \sum_{y\, (\mod{a''})} \e{\frac{b\alpha y}{a''}},
\end{align*}
where
\[ a' = \prod_{v_p(a) = v_p(c)} p^{v_p(c)}, \qquad a''=a/a'. \]
Now the $y$ sum is zero unless $a''=1$, so assume $a=a'$ is coprime to $c/a$, and the sum factors as
\begin{align*}
	S(a n, b\alpha; c) =& \delta_{(a,c/a)=1} S(0, 1; a) S(n, b\alpha\wbar{a}; c/a).
\end{align*}
Again we have
\begin{align*}
	S(n, b\alpha\wbar{a}; c/a) =& \sum_{\substack{x\, (\mod{c/(ab'')})\\(x,c/(ab''))=1}} \e{\frac{n x+b\alpha\wbar{a} \wbar{x}}{c/a}} \sum_{y\, (\mod{b''})} \e{\frac{nx}{b''}},
\end{align*}
where
\[ b' = \prod_{v_p(b) = v_p(c)} p^{v_p(c)}, \qquad b''=b/b', \]
and the sum is zero unless $b''|n$.
Altogether,
\begin{align*}
	S(a n, b\alpha; c) =& \delta_{(a,c/a)=1} \delta_{b''|n} b'' S(0, 1; a) S(0, n; b') S\left( \frac{n}{b''}, \alpha\wbar{ab'}; \frac{c}{ab}\right).
\end{align*}
Applying this to the Kloosterman sum
$$ S\left(E_1 \Bigl(\wbar{\frac{D_1}{e_3}}\Bigr)^2 \frac{D_2}{f_3}, n E_2 \wbar{C_1}; \frac{e_3}{E_3}\right) $$
 in \eqref{multiline}, we see $(E_2, e_3/(E_3 E_2))=1$, but this is impossible unless $$E_2=1$$ since $e_3^\dagger \mid  E_3 \mid \frac{e_3}{E_3 E_2}$ and $E_2 \mid e_3$.
Similarly, we have $$E_4=1.$$
By the same reasoning, we have $$ E_1E_5 \mid (m,n), \quad \text{and} \quad 
 (E_8, e_9/(E_8 E_9))=1.$$

We introduce more variables. Set $J_1 J_2=E_7$, where $$p\mid J_1 \Rightarrow v_p(J_1) = v_p(e_9/E_9), \quad p\mid J_2 \Rightarrow v_p(J_2) < v_p(e_9/E_9). $$  Then in particular  $(J_1, J_2)=1$ and $J_2\mid (m,n)$.
For convenience, set $$J_3 = e_9/(J_1 J_2 E_8 E_9),$$ and substitute $f_j \mapsto f_j e_j$. This gives
\begin{align*}
	& S(1, m, n, 1; D_1, D_2) = S((\wbar{D_1/F_1})^2 D_2, n; F_1) S((\wbar{D_2/F_2})^2 D_1, m; F_2) \\
	& \qquad \times \mathop{{\sum}'}\limits_{\substack{e_2 e_3=G_1\\f_2 f_3 e_2 e_3=G_2\\e_2\mid m}} \sum_{\substack{E_1 E_3^2=e_3 \\ e_3^\dagger \mid E_3 \\ E_1 \mid n, E_1 \mid m}} \sum_{\substack{C_1 \, (\mod{E_3})\\(C_1,E_3)=(C_1^*,f_3 E_3)=1}} e_2 E_1^2 E_3 S(0,n; e_2) S\left(\Bigl(\wbar{D_2/(f_2 e_2)}\Bigr)^2 \frac{D_1}{e_2} \frac{m}{e_2}, 1; f_2\right) \\
	& \qquad\qquad \times S\left(\Bigl(\wbar{\frac{D_1}{e_3}}\Bigr)^2 \frac{D_2}{f_3 e_3}, \frac{n}{E_1} \wbar{C_1}; E_3\right) S\left(\Bigl(\wbar{\frac{D_2}{f_3 e_3}}\Bigr)^2 \frac{D_1}{e_3} \frac{m}{E_1}, \wbar{C_1^*};f_3 E_3\right) \\
	\end{align*}
	\begin{align*}
	& \qquad \times \mathop{{\sum}'}\limits_{\substack{e_4 e_6=H_2\\f_4 f_6 e_4 e_6=H_1\\e_4\mid n}} \sum_{\substack{E_5 E_6^2=e_6 \\ e_6^\dagger \mid E_6 \\ E_5 \mid n, E_5 \mid m}} \sum_{\substack{C_2 \, (\mod{E_6})\\(C_2,E_6)=(C_2^*,f_6 E_6)=1}} e_4 E_5^2 E_6 S(0, m; e_4) S\left(\Bigl(\wbar{\frac{D_1}{f_4 e_4}}\Bigr)^2 \frac{D_2}{e_4}, \frac{n}{e_4}; f_4\right) \\
	& \qquad \qquad\times S\left(\Bigl(\wbar{\frac{D_2}{e_6}}\Bigr)^2 \frac{D_1}{f_6 e_6} \frac{m}{E_5}, \wbar{C_2}; E_6\right) S\left(\Bigl(\wbar{\frac{D_1}{f_6 e_6}}\Bigr)^2 \frac{D_2}{e_6}, \frac{n}{E_5} \wbar{C_2^*};f_6 E_6\right) \\
		& \qquad \times \mathop{{\sum}'}\limits_{\substack{e_7 e_8 e_9=I}} \sum_{\substack{J_1 J_2 J_3 E_8 E_9 = e_9 \\ J_3 \mid E_9, e_9^\dagger \mid E_9 \\ (E_8, J_1 J_2 J_3)=(J_1,J_2 J_3)=1 \\ J_2\mid m, J_2\mid n}} \sum_{\substack{C_3  \, (\mod{J_3})\\(C_3,J_3)=(C_3^*,J_3)=1}} \frac{J_2^2 E_9}{C_{e_9, E_9}} S(0, 1; e_7) S(0,n; e_8)S(0, 1; e_7) S(m, 0; e_8) \\
	& \qquad \times  S(0, 1; E_8)^2 S(0, m; J_1) S(0, n; J_1) S\left(\Bigl(\wbar{\frac{D_1}{e_9}}\Bigr)^2 \frac{D_2}{e_9} \wbar{J_1 E_8}, \frac{n}{J_2} \wbar{C_3}; J_3\right) S\left(\Bigl(\wbar{\frac{D_2}{e_9}}\Bigr)^2 \frac{D_1}{e_9} \frac{m}{J_2}, \wbar{J_1 E_8} \wbar{C_3^*}; J_3\right), 
\end{align*}
where 
\begin{equation}\label{Cstar}
C_1^* = 1-C_1f_3, \quad C_2^* = 1-C_2f_6, \quad C_3^* =  E_9/J_3-C_3.
\end{equation}\\

The key step is now to recognize that one can apply the Chinese Remainder Theorem backwards to re-combine many of the Kloosterman sums. This requires some further  manipulation. The product of the Kloosterman sums to modulus $F_2, f_2, f_3E_3, E_6$ and the second Kloosterman sum to modulus $J_3$ can be written as
\begin{align*}
	& S((m/(e_2 E_1 E_5 J_2)) (\wbar{f_2 f_3 E_3 E_6 J_3})^2, \wbar{e_4 e_7 e_8 J_1 E_8 E_9} (F_1 f_4 f_6 J_3); F_2) \\
	& \qquad \times S((m/(e_2 E_1 E_5 J_2)) (\wbar{F_2 f_3 E_3 E_6 J_3})^2, \wbar{e_4 e_7 e_8 J_1 E_8 E_9} (F_1 f_4 f_6 J_3); f_2) \\
	& \qquad \times S((m/(e_2 E_1 E_5 J_2)) (\wbar{F_2 f_2 E_6 J_3})^2, \wbar{e_4 e_7 e_8 J_1 E_8 E_9} (F_1 f_4 f_6 J_3) \wbar{C_1^*};f_3 E_3) \\
	& \qquad \times S((m/(e_2 E_1 E_5 J_2)) (\wbar{F_2 f_2 f_3 E_3 J_3})^2, \wbar{e_4 e_7 e_8 J_1 E_8 E_9} (F_1 f_4 J_3) \wbar{C_2}; E_6) \\
	& \qquad \times S((m/(e_2 E_1 E_5 J_2)) (\wbar{F_2 f_2 f_3 E_3 E_6})^2, \wbar{e_4 e_7 e_8 J_1 E_8} (F_1 f_4 f_6) \wbar{C_3^*}; J_3) \\
	&= S((m/(e_2 E_1 E_5 J_2)), F_1 \alpha_{F_2}; F_2 f_2 f_3 E_3 E_6 J_3),
\end{align*}
where $\alpha_{F_2}$ is uniquely determined modulo $F_2 f_2 f_3 E_3 E_6 J_3$ by
 \begin{equation}\label{alpha}
\begin{split}
	& \wbar{e_4 e_7 e_8 J_1 E_8 E_9} (f_4 f_6 J_3) \pmod{F_2 f_2}, \\
	& \wbar{e_4 e_7 e_8 J_1 E_8 E_9} (f_4 f_6 J_3) \wbar{C_1^*} \pmod{f_3 E_3}, \\
	& \wbar{e_4 e_7 e_8 J_1 E_8 E_9} (f_4 J_3) \wbar{C_2} \pmod{E_6}, \\
	& \wbar{e_4 e_7 e_8 J_1 E_8} (f_4 f_6) \wbar{C_3^*} \pmod{J_3}.
	\end{split}
	\end{equation}
In particular, $(\alpha_{F_2},F_2 f_2 f_3 E_3 E_6 J_3)=1$ and 
  $\alpha_{F_2}$ depends on $F_2$, but not $F_1$.

Similarly, the product of the Kloosterman sums to modulus $F_1, E_3, f_4, f_6E_6$ and the first instance of the Kloosterman sum to modulus $J_3$ can be written as
\begin{align*}
 	& S((n/(E_1 e_4 E_5 J_2)) (\wbar{E_3 f_4 f_6 E_6 J_3})^2, \wbar{e_2 e_7 e_8 J_1 E_8 E_9} (F_2 f_2 f_3 J_3); F_1) \\
	& \qquad \times S((n/(E_1 e_4 E_5 J_2)) (\wbar{F_1 f_4 f_6 E_6 J_3})^2, \wbar{e_2 e_7 e_8 J_1 E_8 E_9} (F_2 f_2 J_3) \wbar{C_1}; E_3) \\
	& \qquad \times S((n/(E_1 e_4 E_5 J_2)) (\wbar{F_1 E_3 f_6 E_6 J_3})^2, \wbar{e_2 e_7 e_8 J_1 E_8 E_9} (F_2 f_2 f_3 J_3); f_4) \\
	& \qquad \times S((n/(E_1 e_4 E_5 J_2)) (\wbar{F_1 E_3 f_4 J_3})^2, \wbar{e_2 e_7 e_8 J_1 E_8 E_9} (F_2 f_2 f_3 J_3) \wbar{C_2^*};f_6 E_6) \\
	& \qquad \times S((n/(E_1 e_4 E_5 J_2)) (\wbar{F_1 E_3 f_4 f_6 E_6})^2, \wbar{e_2 e_7 e_8 J_1 E_8} (F_2 f_2 f_3) \wbar{C_3}; J_3) \\
	&= S((n/(E_1 e_4 E_5 J_2)), F_2 \beta_{F_1}, F_1 E_3 f_4 f_6 E_6 J_3),
\end{align*}
where $\beta_{F_1}$ is uniquely determined modulo $F_1 E_3 f_4 f_6 E_6 J_3$ by 
\begin{equation}\label{beta}
\begin{split}
	& \wbar{e_2 e_7 e_8 J_1 E_8 E_9} (f_2 f_3 J_3) \pmod{F_1 f_4}, \\
	& \wbar{e_2 e_7 e_8 J_1 E_8 E_9} (f_2 J_3) \wbar{C_1} \pmod{E_3}, \\
	& \wbar{e_2 e_7 e_8 J_1 E_8 E_9} (f_2 f_3 J_3) \wbar{C_2^*} \pmod{f_6 E_6}, \\
	& \wbar{e_2 e_7 e_8 J_1 E_8} (f_2 f_3) \wbar{C_3} \pmod{J_3}.
		\end{split}
	\end{equation}
Again,  $(\beta_{F_1}, F_1 E_3 f_4 f_6 E_6 J_3)=1$  and $\beta_{F_1}$ depends on $F_1$, but not $F_2$. \\

We are now ready to state the complete global decomposition. 
 
\begin{theorem}\label{global} Keep the notation   \eqref{factor},  \eqref{Cstar} and \eqref{prime}.  We have
\begin{align*}
	& S(1, m, n, 1; D_1, D_2) \\
	&= \mathop{{\sum}'}\limits_{\substack{F_1 f_4 f_6 e_2 e_3 e_4 e_6 e_7 e_8 e_9 = D_1 \\ F_2 f_2 f_3 e_2 e_3 e_4 e_6 e_7 e_8 e_9 = D_2 \\ e_2\mid m, e_4\mid n\\ (F_1, F_2)=(F_1 F_2, e_2 e_3 e_4 e_6 e_7 e_8 e_9)=1}} \sum_{\substack{E_1 E_3^2=e_3 \\ e_3^\dagger \mid  E_3 \\ E_1 \mid n, E_1 \mid m}} \sum_{\substack{E_5 E_6^2=e_6 \\ e_6^\dagger | E_6 \\ E_5 \mid n, E_5 \mid m}} \sum_{\substack{J_1 J_2 J_3 E_8 E_9 = e_9 \\ J_3 \mid E_9, e_9^\dagger \mid E_9 \\ (E_8, J_1 J_2 J_3)=(J_1,J_2 J_3)=1 \\ J_2\mid m, J_2\mid n}} \sum_{(\alpha,\beta)\in \mathcal{C}} \frac{e_2 E_1^2 E_3 e_4 E_5^2 E_6 J_2^2 E_9}{C_{e_9, E_9}} \\
	& \qquad \times \mu(e_7 E_8)^2 S(0, m; e_4 e_8 J_1) S(0, n; e_2 e_8 J_1) \\
	& \qquad \times S\left(\frac{m}{e_2 E_1 E_5 J_2}, F_1 \alpha_{F_2}; F_2 f_2 f_3 E_3 E_6 J_3\right) S\left(\frac{n}{E_1 e_4 E_5 J_2}, F_2 \beta_{F_1}, F_1 E_3 f_4 f_6 E_6 J_3\right), 
\end{align*}
where $\mu$ is the M\"obius function, 
$\mathcal{C} = \mathcal{C}(e_2, f_2, f_3, e_4, f_4, f_6, e_7, e_8, E_3, E_6, E_7, E_8, E_9, J_1, J_2, J_3)$
is in bijection with triples
\begin{align*}
	&\set{C_1 \, (\mod{E_3})\mid (C_1,E_3)=(C_1^*,f_3 E_3)=1} \\\
	&  \times \set{C_2 \, (\mod{E_6})  \mid (C_2,E_6)=(C_2^*,f_6 E_6)=1} \\
	&   \times \set{C_3 \, (\mod{J_3}) \mid (C_3,J_3)=(C_3^*,J_3)=1} , 
\end{align*}
and each $(C_1, C_2, C_3)$ determines $(\alpha_{F_2}, \beta_{F_1})$ via 
  \eqref{alpha} and \eqref{beta}.  \end{theorem} 
 Of course, the same type of decomposition holds for $S(\pm 1, \pm m, n, 1; D_1, D_2).$
 
 This is the announced generalization of \eqref{klooster1} and \eqref{klooster2} (for $n_1 = n_2 = 1$). For instance, if $(D_1, D_2) = 1$, then $D_1 = F_1$, $D_2 = F_2$, all other variables are 1 (in particular $\mathcal{C}$ is trivial), and $\alpha_{F_2} =  \beta_{F_1} = 1$. On the other hand, if $D_1 = D_2 = p$ is a prime not dividing $nm$, then the only possibilities are $e_7 = p$ and all other variables are 1, $e_8 = p$ and all other variables are 1, or $e_9 = E_9 = p$ and all other variables are 1, and we obtain
 $$S(1, m, n, 1; p, p) = 1 + 1 + \frac{p}{C_{p, p}} = p+1.$$
    
Although admittedly rather complicated, we will see in the next section how to make good use of this formula. Roughly speaking it makes rigorous the heuristic idea that reality is not much harder than \eqref{klooster1} and \eqref{klooster2} together. 
    
\section{Bilinear forms in Kloosterman sums}

In this section we prove Theorem \ref{thm1}. We keep the notation from the previous section. For later purposes we consider a slightly more general quantity than the sum considered in Theorem \ref{thm1}. For $s = (s_1, s_2) $ with $\Re s_1 = \Re s_2 = 0$ let $\Phi(s)$ be a  function  with support in $|s_1| \leq T_1$, $|s_2| \leq T_2$, and suppose that
$$\sup_{|s_2| \leq T_2} \int_{\Re s_1 = 0} |\Phi(s)| ds_1 \leq S_2, \quad \sup_{|s_1| \leq T_1} \int_{\Re s_2 = 0} |\Phi(s)| ds_2 \leq S_1.$$
Then we consider
$$\mathcal{S} = \int_{\Re s = 0} |\Phi(s)|  \sum_{\substack{D_1 \le X_1\\ D_2 \le X_2}} \Bigl|\sum_{m,n\le N} a_nb_m n^{-s_1} m^{-s_2} S(1,m,n,1;D_1,D_2)\Bigr| ds $$
with the aim of proving
\begin{prop}\label{hybrid} With the above notation we have
$$\frac{\mathcal{S}}{X_1X_2}  \ll (S_1S_2)^{1/2} (X_1X_2)^{\varepsilon} \| a \|_2 \| b \|_2 \Biggl((X_1 + X_2) (T_1T_2)^{1/2} + \frac{X_1 T_1^{1/2} + X_2 T_2^{1/2}}{\min(X_1, X_2)^{1/2}}\sqrt{N} +  N\Biggr).$$
\end{prop}

\textbf{Proof.} By the Cauchy-Schwarz inequality and Theorem \ref{global} we have
\begin{equation}\label{S}
\begin{split}
\mathcal{S}	\le& \mathop{{\sum}'}\limits_{\substack{f_4 f_6 e_2 e_3 e_4 e_6 e_7 e_8 e_9 \le X_1 \\ f_2 f_3 e_2 e_3 e_4 e_6 e_7 e_8 e_9 \le X_2}} \sum_{\substack{E_1 E_3^2=e_3 \\ e_3^\dagger \mid E_3 \\ E_1 \mid n, E_1 \mid m}} \sum_{\substack{E_5 E_6^2=e_6 \\ e_6^\dagger \mid E_6 \\ E_5 \mid n, E_5 \mid m}} \sum_{\substack{J_1 J_2 J_3 E_8 E_9 = e_9 \\ J_3 \mid E_9, e_9^\dagger \mid E_9 \\ (E_8, J_1 J_2 J_3)=(J_1,J_2 J_3)=1 \\ J_2\mid m, J_2\mid n}} \sum_{(\alpha,\beta)\in \mathcal{C}} \\
	& \qquad \times \sum_{\substack{g_1\mid e_4 e_8 J_1\\g_2\mid e_2 e_8 J_1}} g_1 g_2 e_2 E_1^2 E_3 e_4 E_5^2 E_6 J_2^2 E_9 \sqrt{U_1 U_2},
 \end{split}
 \end{equation}
 where
\begin{align*}
	U_1 =& \int_{\Re s = 0} |\Phi(s)| \sum_{\substack{F_1 \le X_1/(f_4 f_6 e_2 e_3 e_4 e_6 e_7 e_8 e_9) \\ F_2 \le X_2/(f_2 f_3 e_2 e_3 e_4 e_6 e_7 e_8 e_9) \\ (F_1, F_2)=(F_1 F_2, e_2 e_3 e_4 e_6 e_7 e_8 e_9)=1}} \\
	& \quad\quad\quad\quad\times \Bigl|\sum_{\substack{m\le N/(e_2 E_1 E_5 J_2)\\(m,e_4 e_8 J_1)=g_1}} b_{e_2 E_1 E_5 J_2 m}m^{-s_2} S(m, F_1 \alpha_{F_2}; F_2 f_2 f_3 E_3 E_6 J_3)\Bigr|^2 ds,
\end{align*}
\begin{align*}
	U_2 =&\int_{\Re s = 0} |\Phi(s)|  \sum_{\substack{F_1 \le X_1/(f_4 f_6 e_2 e_3 e_4 e_6 e_7 e_8 e_9) \\ F_2 \le X_2/(f_2 f_3 e_2 e_3 e_4 e_6 e_7 e_8 e_9) \\ (F_1, F_2)=(F_1 F_2, e_2 e_3 e_4 e_6 e_7 e_8 e_9)=1}}\\
	&\quad\quad\quad\quad\times  \Bigl|\sum_{\substack{n\le N/(e_4 E_1 E_5 J_2)\\(n,e_2 e_8 J_1)=g_2}} a_{e_4 E_1 E_5 J_2 n} n^{-s_1} S(n, F_2 \beta_{F_1}, F_1 E_3 f_4 f_6 E_6 J_3)\Bigr|^2 ds.
\end{align*}
In the expression for $U_1$, we drop the condition $(F_1, F_2 e_2 e_3 e_4 e_6 e_7 e_8 e_9)=1$, let $$k=\ceil{\frac{X_1/(f_4 f_6 e_2 e_3 e_4 e_6 e_7 e_8 e_9)}{F_2 f_2 f_3 E_3 E_6 J_3}}$$ and extend the $F_1$ sum to $k F_2 f_2 f_3 E_3 E_6 J_3$. Keeping in mind that $F_2 f_2 f_3 E_3 E_6 J_3 \leq X_2  E_3 E_6 J_3/(e_2 e_3 e_4 e_6 e_7 e_8 e_9)$ and $e_3 = E_1E_3^2$, $e_6 = E_5E_6^2$, $e_9 = J_1J_2J_3 E_8E_9$, we see that  the $F_1$ sum is at most
\begin{align*}
	& \sum_{F_1 \le k F_2 f_2 f_3 E_3 E_6 J_3} \Biggl|\sum_{\substack{m\le N/(e_2 E_1 E_5 J_2)\\(m,e_4 e_8 J_1)=g_1}}b_{e_2 E_1 E_5 J_2 m} m^{-s_2} S(m, F_1 \alpha_{F_2}; F_2 f_2 f_3 E_3 E_6 J_3)\Biggr|^2 \\
		\end{align*}
	\begin{align*}
	&= \sum_{\substack{m_1, m_2\le N/(e_2 E_1 E_5 J_2)\\(m_1,e_4 e_8 J_1)=(m_2,e_4 e_8 J_1)=g_1}} b_{e_2 E_1 E_5 J_2 m_1} \wbar{b_{e_2 E_1 E_5 J_2 m_2}} \left(\frac{m_1}{m_2}\right)^{-s_2}  \\
	& \quad\quad  \sum_{\substack{x_1, x_2 \, (\mod{F_2 f_2 f_3 E_3 E_6 J_3})\\(x_1 x_2, F_2 f_2 f_3 E_3 E_6 J_3)=1}} \e{\frac{m_1 x_1-m_2 x_2}{F_2 f_2 f_3 E_3 E_6 J_3}} \sum_{F_1 \le k F_2 f_2 f_3 E_3 E_6 J_3} \e{\frac{F_1 \alpha_{F_2}(\wbar{x_1}-\wbar{x_2})}{F_2 f_2 f_3 E_3 E_6 J_3}} \\
	&\le \left(\frac{X_1}{f_4 f_6 e_2 e_3 e_4 e_6 e_7 e_8 e_9}+\frac{X_2}{e_2 E_1 E_3 e_4 E_5 E_6 e_7 e_8 J_1 J_2 E_8 E_9}\right) \\
	&    \sum_{\substack{m_1, m_2\le N/(e_2 E_1 E_5 J_2)\\(m_1,e_4 e_8 J_1)=(m_2,e_4 e_8 J_1)=g_1}} b_{e_2 E_1 E_5 J_2 m_1} \wbar{b_{e_2 E_1 E_5 J_2 m_2}} \left(\frac{m_1}{m_2}\right)^{-s_2}  \sum_{\substack{x \, (\mod{F_2 f_2 f_3 E_3 E_6 J_3})\\(x, F_2 f_2 f_3 E_3 E_6 J_3)=1}} \e{\frac{(m_1-m_2) x}{F_2 f_2 f_3 E_3 E_6 J_3}} \\
		&= \left(\frac{X_1}{f_4 f_6 e_2 e_3 e_4 e_6 e_7 e_8 e_9}+\frac{X_2}{e_2 E_1 E_3 e_4 E_5 E_6 e_7 e_8 J_1 J_2 E_8 E_9}\right) \\
	&     \sum_{\substack{x \, (\mod{F_2 f_2 f_3 E_3 E_6 J_3})\\(x, F_2 f_2 f_3 E_3 E_6 J_3)=1}} \Biggl|\sum_{\substack{m\le N/(e_2 E_1 E_5 J_2)\\(m,e_4 e_8 J_1)=g_1}} b_{e_2 E_1 E_5 J_2 m} m^{-s_2} \e{\frac{m x}{F_2 f_2 f_3 E_3 E_6 J_3}}\Biggr|^2.
\end{align*}
Using also that $(g_1, F_2 f_2 f_3 E_3 E_6 J_3)=1$, we obtain
\begin{align*}
	U_1 \le& S_1 \left(\frac{X_1}{f_4 f_6 e_2 e_3 e_4 e_6 e_7 e_8 e_9}+\frac{X_2}{e_2 E_1 E_3 e_4 E_5 E_6 e_7 e_8 J_1 J_2 E_8 E_9}\right)\sum_{\substack{F_2 \le X_2/(f_2 f_3 e_2 e_3 e_4 e_6 e_7 e_8 e_9) \\ (F_2, e_2 e_3 e_4 e_6 e_7 e_8 e_9)=1}}  \\
	&    \sum_{\substack{x \, (\mod{F_2 f_2 f_3 E_3 E_6 J_3})\\(x, F_2 f_2 f_3 E_3 E_6 J_3)=1}} \int_{-T_2}^{T_2}  \Biggl|\sum_{\substack{m\le N/(e_2 E_1 E_5 J_2g_1)\\(m,e_4 e_8 J_1/g_1)=1}} b_{e_2 E_1 E_5 J_2 g_1 m}  m^{-it} \e{\frac{m x}{F_2 f_2 f_3 E_3 E_6 J_3}}\Biggr|^2 dt. 
\end{align*}
We now invoke a slight variation of Gallagher's hybrid large sieve inequality \cite{Ga} whose proof is almost verbatim the same (just use the additive large sieve \cite[Theorem 7.7]{IK} instead of the multiplicative version \cite[Lemma 3]{Ga}):

\begin{lemma}\label{gallagher} For $f \in \Bbb{N}$, $X , N, T \geq 1$ and $b(n)$, $n \leq N$, any sequence of complex numbers, we have
$$\sum_{F \leq X} \sum_{x \, (\text{{\rm mod }} Ff)} \int_{-T}^T \Bigl| \sum_{m \leq M} b(m) m^{-it} e\left(\frac{mx}{Ff}\right) \Bigr|^2 dt \ll (X^2 f T + N)\| b \|_2^2. $$
\end{lemma} 

This gives us
\begin{align*}
	U_1 \ll& S_1 \left(\frac{X_1}{f_4 f_6 e_2 e_3 e_4 e_6 e_7 e_8 e_9}+\frac{X_2}{e_2 E_1 E_3 e_4 E_5 E_6 e_7 e_8 J_1 J_2 E_8 E_9}\right) \\
	& \qquad \times \left(\frac{X_2^2T_2}{f_2 f_3 e_2^2 E_1^2 E_3^3 e_4^2 E_5^2 E_6^3 e_7^2 e_8^2 J_1^2 J_2^2 J_3 E_8^2 E_9^2} + \frac{N}{e_2 E_1 E_5 J_2 g_1}\right) \norm{b}_2^2
\end{align*}
and analogously 
\begin{align*}
	U_2 \ll& S_2 \left(\frac{X_2}{f_2 f_3 e_2 e_3 e_4 e_6 e_7 e_8 e_9}+\frac{X_1}{e_2 E_1 E_3 e_4 E_5 E_6 e_7 e_8 J_1 J_2 E_8 E_9}\right) \\
	& \qquad \times \left(\frac{X_1^2T_1}{f_4 f_6 e_2^2 E_1^2 E_3^3 e_4^2 E_5^2 E_6^3 e_7^2 e_8^2 J_1^2 J_2^2 J_3 E_8^2 E_9^2} + \frac{N}{e_4 E_1 E_5 J_2 g_2}\right) \norm{a}_2^2.
\end{align*}
It is now a matter of book-keeping. Substituting back into \eqref{S} and using $\left|\mathcal{C}\right| \le E_3 E_6 J_3$, we obtain
\begin{align*}
	\mathcal{S} \ll& (S_1S_2)^{1/2} \norm{a}_2 \norm{b}_2 \mathop{{\sum}'}\limits_{\substack{f_4 f_6 e_2 e_3 e_4 e_6 e_7 e_8 e_9 \le X_1 \\ f_2 f_3 e_2 e_3 e_4 e_6 e_7 e_8 e_9 \le X_2}} \sum_{\substack{E_1, E_3, E_5, E_6 \\ J_1, J_2, J_3, E_8, E_9}} \sum_{\substack{g_1|e_4 e_8 J_1\\g_2|e_2 e_8 J_1}} \frac{g_1 g_2 E_3 E_6 J_3}{\sqrt{e_2 e_4} e_7 e_8 J_1 E_8} \\
	& \qquad \times \left(\sqrt{X_1 X_2}+\frac{X_1+X_2}{\sqrt{f_2 f_3 f_4 f_6} E_3 E_6 J_3}\right)  \left(\frac{X_2^2T_2}{f_2 f_3 e_2 E_1 E_3^3 e_4^2 E_5 E_6^3 e_7^2 e_8^2 J_1^2 J_2 J_3 E_8^2 E_9^2} + \frac{N}{g_1}\right)^{1/2}\\
	& \qquad  \times \left(\frac{X_1^2T_1}{f_4 f_6 e_2^2 E_1 E_3^3 e_4 E_5 E_6^3 e_7^2 e_8^2 J_1^2 J_2 J_3 E_8^2 E_9^2} + \frac{N}{g_2}\right)^{1/2},
\end{align*}
where $E_1, E_3, E_5, E_6, E_8, E_9, J_1, J_2, J_3$ are subject to the same conditions as in \eqref{S}. Given $e_2$, we see by Rankin's trick that
$$\underset{f_2 \leq X}{\left.\sum\right.'} 1 \leq X^{\varepsilon} \left.\sum_{f_2}\right.' \frac{1}{f_2^{\varepsilon}} = X^{\varepsilon} \prod_{p \mid e_2}  \frac{p^{\varepsilon}}{p^{\varepsilon} - 1} \ll (Xe_2)^{\varepsilon}.$$
Hence summing over $f_2, f_3, f_4, f_6, e_7, g_1, g_2$, we obtain
\begin{align*}
	\mathcal{S} \ll& (S_1S_2)^{1/2} (X_1 X_2)^\varepsilon \norm{a}_2 \norm{b}_2 \mathop{{\sum}'}\limits_{e_2 e_3 e_4 e_6 e_8 e_9 \le \min(X_1, X_2)} \sum_{\substack{E_1, E_3, E_5, E_6 \\ J_1, J_2, J_3, E_8, E_9}} \frac{E_3 E_6 J_3}{E_8} \\
	& \qquad \times \left(\sqrt{X_1 X_2}+\frac{X_1+X_2}{E_3 E_6 J_3}\right) \left(\frac{X_2T_2^{1/2}}{(e_2  E_1  E_3^3 e_4 E_5  E_6^3 e_8  J_1  J_2  J_3)^{1/2}  E_8 E_9} + \sqrt{N}\right) \\
	& \qquad \times \left(\frac{X_1T_1^{1/2}}{(e_2 E_1  E_3^{3} e_4  E_5  E_6^{3} e_8  J_1 J_2 J_3)^{1/2} E_8 E_9} + \sqrt{N}\right).
\end{align*}
Next we sum over $E_1, E_5, E_8, J_1, J_2$ getting
\begin{align*}
	\mathcal{S} \ll& (S_1S_2)^{1/2} (X_1 X_2)^\varepsilon \norm{a}_2 \norm{b}_2 \sum_{\substack{e_2 E_3^2 e_4 E_5^2 e_8 J_3 E_9 \le \min(X_1, X_2)\\J_3 | E_9}} E_3 E_6 J_3 \left(\sqrt{X_1 X_2}+\frac{X_1+X_2}{E_3 E_6 J_3}\right) \\
	& \qquad \times \left(\frac{X_1 X_2 (T_1T_2)^{1/2}}{e_2 E_3^3 e_4 E_6^3 e_8 J_3 E_9^2} + \frac{(X_1T_1^{1/2}+X_2T_2^{1/2})\sqrt{N}}{(e_2  E_3^{3 } e_4  E_6^{3 } e_8  J_3)^{1/2} E_9} + N\right).
\end{align*}
Execute the $e_2$, $e_4$, $e_8$ sums
\begin{align*}
	\mathcal{S} \ll& (S_1S_2)^{1/2} (X_1 X_2)^\varepsilon \norm{a}_2 \norm{b}_2 \sum_{\substack{E_3^2 E_6^2 J_3 E_9 \le \min(X_1, X_2)\\J_3 | E_9}} E_3 E_6 J_3 \left(\sqrt{X_1 X_2}+\frac{X_1+X_2}{E_3 E_6 J_3}\right) \\
	& \qquad \times \left(\frac{X_1 X_2(T_1T_2)^{1/2} }{E_3^3 E_6^3 J_3 E_9^2} + \frac{(X_1T_1^{1/2}+X_2T_2^{1/2}) \min(X_1, X_2)^{1/2}  N}{E_3^{5/2} E_6^{5/2} J_3 E_9^{3/2}} + \frac{N \min(X_1, X_2)}{E_3^2 E_6^2 J_3 E_9}\right).
\end{align*}
Finally summing over $E_3, E_6, J_3$ and $E_9$ completes the proof of Proposition \ref{hybrid}. \hfill $\square$ \\

The proof of   Theorem \ref{thm1} is almost verbatim the same. We simply ignore the integration over $s$, put $T_1 = T_2 = S_1 = S_2 = 1$ and use the non-hybrid large sieve inequality (\cite[Theorem 7.7]{IK})
$$\sum_{F \leq X} \sum_{x \, (\text{{\rm mod }} Ff)}  \Bigl| \sum_{m \leq M} b(m)  e\left(\frac{mx}{Ff}\right) \Bigr|^2 dt \leq (X^2 f + N)\| b \|_2^2 $$
 instead of Lemma \ref{gallagher}. \hfill $\square$

\section{Proof of Proposition \ref{lowerbound}}\label{proof}

The lower bound in Proposition \ref{lowerbound}  comes from the maximal Eisenstein series $ E(z, 1/2 + it; u_j)$. These are parametrized by $t \in \Bbb{R}$ and   Hecke-Maa{\ss} cusp forms $u_j$ for the group ${\rm SL}_2(\Bbb{Z})$ with spectral parameter $t_j$ and Hecke eigenvalues $\lambda_j(n)$. See \cite[Section 10]{Go} for details. 
 The Hecke eigenvalues of $E(z, 1/2 + it; u_j)$ are given by 
$$\lambda(n)  = \sum_{d_1d_2 = n} \lambda_{j}(d_1) d_1^{-it} d_2^{2 it}$$
and the spectral parameter is $\mu = \mu(t, t_j) = (2it, -it + it_j, -it-it_j)$. 
The normalizing factor $\mathcal{N}(\pi)$ is proportional to $L(1, \text{Ad}^2 u_j) |L(1+ 3 i t, u_j)|^2$ (see \cite[Section 3.1]{BB}). 

Choose some ball $B \subseteq \Omega$ and let $\tau \in \Bbb{R}$ be any real number such that $\mu(\tau, t) \in B$ for some $t \in \Bbb{R}$ (and hence for a small interval of $t$). Then we choose $a_n = W(n/N) (n/N)^{2iT\tau}$ for a fixed,  non-negative, non-zero smooth  weight function $W$ with support in $[1, 2]$. Hence the maximal Eisenstein contribution is of the shape 
$$\underset{\mu(t, t_j) \in T\Omega} {\int   \sum} \frac{1}{L(1, \text{Ad}^2 u_j) |L(1+ 3 i t, u_j)|^2}   \Bigl| \sum_{d_1, d_2} \lambda_j(d_1) d_1^{-it} d_2^{2 it}  W\Bigl(\frac{d_1d_2}{N}\Bigr)\Bigl(\frac{d_1d_2}{N}\Bigr)^{2iT\tau}  \Bigr|^2 dt.$$
By Mellin inversion, the inner sum equals
$$\frac{1}{2\pi i} \int_{(2)} \widehat{W}(s-2iT\tau) N^s \zeta(s - 2it) L( s + it, u_j) ds.$$
The analytic conductor of $\zeta(s - 2it)L(u_j, s + it)$ is $O(T^3)$  for $\mu(t, t_j) \in T\Omega$. Since  $N \geq T^{3 + \delta}$, we can shift contours to the left as far as we wish and pick up a pole at $s = 1+ 2 it$ whose residual contribution is
$$\underset{\mu(t, t_j) \in T\Omega} {\int   \sum} \frac{1}{L(1, \text{Ad}^2 u_j)}  \Bigl| \widehat{W}(1 + 2it - 2iT\tau ) N^{1 + 2 it}  \Bigr|^2 dt \gg N^2   \sum_{\mu(T\tau, t_j) \in T\Omega}   \frac{1}{L(1, \text{Ad}^2 u_j)}. $$
By Weyl's law for ${\rm GL}(2)$, more precisely by the ${\rm GL}(2)$ Kuznetsov formula (\cite[Theorem 16.8]{IK}) to account for the weight function $ L(1, \text{Ad}^2 u_j)^{-1}$, this is 
$$\gg N^2 T^2 \asymp T^2 N \| a \|_2^2,$$
as claimed. \hfill $\square$

\section{The Kuznetsov formula}\label{kuznetsov}

\subsection{Statement of the formula} In this section we state the  Kuznetsov formula.  
We will be brief and refer to \cite[Section 3]{BB} for more details and notation.  
It is convenient to also use coordinates $\nu = (\nu_1, \nu_2, \nu_3) \in i\mathfrak{a}^{\ast}$ given by
 \begin{equation*}\label{nu}
 \nu_1 = \frac{1}{3}(\mu_1 - \mu_2), \quad \nu_2 = \frac{1}{3}(\mu_2 - \mu_3), \quad \nu_3 = \frac{1}{3}(\mu_3 - \mu_1).
 \end{equation*}
The long Weyl element Kloosterman sum was defined in \eqref{klooster0}. In addition, we need another type of Kloosterman given by
\begin{displaymath}
	\tilde{S}(n_1,n_2,m_1;D_1,D_2) := \sum_{\substack{C_1 (\text{mod }D_1), C_2 (\text{mod }D_2)\\(C_1,D_1)=(C_2,D_2/D_1)=1}} e\left(n_2\frac{\bar{C_1}C_2}{D_1}+m_1\frac{\bar{C_2}}{D_2/D_1}+n_1\frac{C_1}{D_1}\right)
\end{displaymath}
for $D_1\mid D_2$. Next,  
let $\mathcal{W}$ 
be the Weyl group. Let 
$h$ be a function  that is holomorphic on $$\{\mu \in \Bbb{C}^3 \mid \mu_1 + \mu_2 + \mu_3 = 0,\, |\Re \mu_j| \leq 3/4\},$$ symmetric under the Weyl group, rapidly decaying as $|\Im \mu_j| \rightarrow \infty$  and satisfies 
\begin{equation*}\label{zeros}
  h(\mu) = 0 \quad \text{whenever} \quad   3\nu_j = \pm 1, \quad j = 1, 2, 3.
\end{equation*}
 We define the spectral measure by
\begin{equation*}
\text{spec}(\mu) d\mu,  \quad  \text{spec}(\mu) := \prod_{j=1}^3 \left(3\nu_j \tan\Bigl(\frac{3\pi}{2} \nu_j\Bigr)\right),
\end{equation*}
where $d\mu = d\mu_1 d\mu_2= d\mu_1d\mu_3 = d\mu_2 d\mu_3$ is the standard measure on the hyperplane $\mu_1 + \mu_2 + \mu_3 = 0$. The Fourier coefficients of a (not necessarily cuspidal) automorphic form $\phi$ generating a representation $\pi$ are (arithmetically) normalized   by $ \lambda_{\pi}(m) = A_{\pi}(1, m) = \overline{A_{\pi}(m, 1)}$.  
Then for $n_1, n_2, m_1, m_2 \in \Bbb{N}$ and $h$ as above we have 
\begin{displaymath}
\begin{split}
&  \int \overline{A_{\pi}(m_1, m_2)}A_{\pi}(n_1, n_2) \frac{h(\mu_{\pi})}{\mathcal{N}(\pi)}  d\pi = \Delta + \Sigma_4 + \Sigma_5 + \Sigma_6,
\end{split}
\end{displaymath}
with
\begin{displaymath}
\begin{split}
  \Delta& = \delta_{n_1, m_1} \delta_{n_2, m_2}  \frac{1}{192\pi^5} \int_{\Re \mu = 0} h(\mu) \text{spec}(\mu) d \mu,\\
   \Sigma_{4}& = \sum_{\epsilon  = \pm 1} \sum_{\substack{D_2 \mid D_1\\  m_2 D_1= n_1 D_2^2}}\frac{ \tilde{S}(-\epsilon n_2, m_2, m_1; D_2, D_1)}{D_1D_2} \Phi_{w_4}\left(  \frac{\epsilon m_1m_2n_2}{D_1 D_2} \right),  \\ 
 \Sigma_{5} &= \sum_{\epsilon  = \pm 1} \sum_{\substack{     D_1 \mid D_2\\ m_1 D_2 = n_2 D_1^2}} \frac{ \tilde{S}(\epsilon n_1, m_1, m_2; D_1, D_2) }{D_1D_2}\Phi_{w_5}\left( \frac{\epsilon n_1m_1m_2}{D_1 D_2}\right),\\
\end{split}
\end{displaymath}
and 
\begin{equation}\label{longWeyl}
   \Sigma_6 = \sum_{\epsilon_1, \epsilon_2 = \pm 1} \sum_{D_1,  D_2  } \frac{S(\epsilon_2 n_2, \epsilon_1 n_1, m_1, m_2; D_1, D_2)}{D_1D_2} \Phi_{w_6}  \left( - \frac{\epsilon_2 m_1n_2D_2}{D_1^2}, - \frac{\epsilon_1 m_2n_1D_1}{ D_2^2}\right),
 \end{equation}
where
\begin{equation}\label{defPhi}
\begin{split}
& \Phi_{w_4}(y) =  \int_{\Re \mu = 0} h(\mu) K_{w_4}(y; \mu )\, \text{spec}(\mu) d \mu,\\
& \Phi_{w_5}(y) = \int_{\Re \mu = 0} h(\mu) K_{w_4}(-y; -\mu )\, \text{spec}(\mu) d \mu ,\\
& \Phi_{w_6}(y_1, y_2) = \int_{\Re \mu = 0} h(\mu) K^{\text{sgn}(y_1), \text{sgn}(y_2)}_{w_6}((y_1, y_2) ;  \mu )\, \text{spec}(\mu) d \mu
\end{split}
\end{equation}
for certain kernel functions $K_{w_4}$, $K^{\pm, \pm}_{w_6}$ whose properties we are going to describe in a moment. 

\subsection{Choice of test function} We will choose the following test function $h$, which approximates the characteristic function on $T\Omega$.  Let $\mu_0 \in \Omega$. We put   $\psi(\mu) = \exp\left(\mu_1^2 +\mu_2^2  + \mu_3^2\right)$ and as in \cite[Section 3.5]{BB} we put 
$$P(\mu) :=  \prod_{0 \leq n \leq A} \prod_{j=1}^3 \frac{(\nu_j - \frac{1}{3}(1 + 2n))(\nu_j + \frac{1}{3}(1 + 2n))}{|\nu_{0, j}|^2} $$
for some large, fixed constant $A$ to compensate  poles of the spectral measure in a large tube.  
Now we  choose
\begin{equation}\label{defh}
h(\mu) :=  P(\mu)^2 \Bigl(\sum_{w \in \mathcal{W}}\psi\Bigl(\frac{w(\mu)  -  T\mu_0}{T^{1-\varepsilon}}\Bigr)\Bigr)^2
\end{equation}
for some very small $\varepsilon$. Then $T^{\varepsilon}$ of these functions are a majorant of $\textbf{1}_{T\Omega}$.

The holomorphicity of $h$ and the location of its zeros are necessary to make the arithmetic side of the Kuznetsov formula absolutely convergent, and in particular to truncate the $D_1, D_2$-sum in \eqref{longWeyl} at the cost of  a negligible error. Having done this, it is convenient to replace $h$, up to a negligible error,   with a real-analytic, Weyl group invariant  function $\tilde{h}$ that is \emph{compactly supported} in $\Omega' T$, where $\Omega' \supseteq \Omega$ is a slightly bigger compact subset not intersecting the Weyl chamber walls, and  $\tilde{h}$ satisfies
\begin{equation}\label{diffh}
\mathscr{D}_j \tilde{h} \ll T^{j(\varepsilon-1)}
\end{equation}
 for every differential operator of order $j$.  We will use this bound frequently when we integrate by parts. Moreover,
\begin{equation}\label{volume}
\int_{\Re \mu = 0} \tilde{h}(\mu) \text{spec}(\mu) d\mu \ll T^5.
\end{equation}

We have already seen in \cite[Lemma 9]{BB} that even a  spectral average of $K_{w_6}^{\pm, \pm}(y; \mu)$ over a $T^{\varepsilon}$-ball in the $\mu$-plane makes $\Phi_{w_6}(y)$ negligibly small 
unless 
\begin{equation}\label{unless}
  \min(|y_2|^{1/3} |y_1|^{1/6},  |y_1|^{1/3} |y_2|^{1/6}) \gg T^{1-\varepsilon}. 
\end{equation}
This   proof relied only on the fact that the center $\mu_0$ of the ball is away from the walls of Weyl chambers (i.e.\   $|\nu_j| \asymp T$ for $j = 1, 2, 3$). In particular, it holds a fortiori for our present test function \eqref{defh}, and so we conclude in the present case  $\Phi_{j}(y) \ll T^{-B}$ unless \eqref{unless} holds.  


\subsection{Kernel functions}\label{kernel} For $x > 0$, $\alpha \in \Bbb{C}$ let
\begin{equation*}\label{jplus}
J_{\alpha}^+(x) := \frac{\pi}{2} \frac{J_{-\alpha}(2x) + J_{\alpha}(2x)}{\cos(\pi \alpha/2)}, \quad J_{\alpha}^-(x) := \frac{\pi}{2} \frac{J_{-\alpha}(2x) - J_{\alpha}(2x)}{\sin(\pi \alpha/2)}, \quad \tilde{K}_{\alpha}(x) = 2\cos\left(\frac{\pi}{2} \alpha\right) K_{\alpha}(2x),
\end{equation*}
where $J_{\alpha}$ and $K_{\alpha}$ are the usual Bessel functions. For $s = (s_1, s_2) \in \Bbb{C}^2$, $\mu \in \Bbb{C}^3$ with $\mu_1 + \mu_2 + \mu_3 = 0$  define the meromorphic function
\begin{equation}\label{Gsmu}
G(s, \mu) :=  \frac{1}{\Gamma(s_1 + s_2)} \prod_{j=1}^3 \Gamma(s_1 - \mu_j) \Gamma(s_2 + \mu_j)\end{equation}
and  the following trigonometric functions
\begin{displaymath}
\begin{split}
& S^{++}(s, \mu) := \frac{1}{24 \pi^2} \prod_{j=1}^3 \cos\left(\frac{3}{2} \pi \nu_j\right),\\
&  S^{+-}(s, \mu) :=  -\frac{1}{32 \pi^2} \frac{\cos(\frac{3}{2} \pi \nu_2)\sin(\pi(s_1 - \mu_1))\sin(\pi(s_2 + \mu_2))\sin(\pi(s_2 + \mu_3))}{\sin(\frac{3}{2} \pi \nu_1)\sin(\frac{3}{2} \pi \nu_3) \sin(\pi(s_1+s_2))}, \\
& S^{-+}(s, \mu) :=-\frac{1}{32 \pi^2}  \frac{\cos(\frac{3}{2} \pi \nu_1)\sin(\pi(s_1 - \mu_1))\sin(\pi(s_1 - \mu_2))\sin(\pi(s_2 + \mu_3))}{\sin(\frac{3}{2} \pi \nu_2)\sin(\frac{3}{2} \pi \nu_3)\sin(\pi(s_1+s_2))}, \\
& S^{--}(s, \mu) := \frac{1}{32 \pi^2}  \frac{\cos(\frac{3}{2} \pi \nu_3) \sin(\pi(s_1 - \mu_2))\sin(\pi(s_2 + \mu_2))}{\sin(\frac{3}{2} \pi \nu_2)\sin(\frac{3}{2} \pi \nu_1)}. 
  \end{split}
\end{displaymath}
Then for   $y = (y_1, y_2) \in (\Bbb{R}\setminus \{0\})^2$ with $\text{{\rm sgn}}(y_1) = \epsilon_1$, $\text{{\rm sgn}}(y_2) = \epsilon_2$ the kernel function $K^{\epsilon_1, \epsilon_2}_{w_6}(y; \mu)$ in \eqref{defPhi} is given by
\begin{equation}\label{defK}
\begin{split}
 K^{\epsilon_1, \epsilon_2}_{w_6}(y; \mu)  =    & \int_{-i\infty}^{i\infty} \int_{-i\infty}^{i\infty}  |4\pi^2 y_1|^{-s_1} |4\pi^2 y_2|^{-s_2}  G(s, \mu) S^{\epsilon_1, \epsilon_2}(s, \mu)\frac{ds_1\, ds_2}{(2\pi i)^2}. 
\end{split}
\end{equation}
The path of integration has to be chosen according to the Barnes convention as in \cite[Definition 1]{BB}. 
As proved in \cite{BB}, there is an alternative description in terms of double Bessel integrals. We define
\begin{equation}\label{BesselJ1}
 \mathcal{J}^{\pm}_{1}(y; \mu) = 
  \Bigl| \frac{y_1}{y_2}\Bigr|^{\frac{1}{2}\mu_2}      \int_0^\infty J^{\pm}_{3\nu_3 }\left(2\pi |y_1|^{1/2}\sqrt{1+u^2}\right) J^{\pm}_{3\nu_3}\left(2\pi |y_2|^{1/2}\sqrt{1+u^{-2}}\right) u^{3\mu_2} \frac{du}{u},
 \end{equation}
\begin{equation}\label{BesselJ2}
 \mathcal{J}_{2}(y; \mu) = 
 \Bigl| \frac{y_1}{y_2}\Bigr|^{\frac{1}{2}\mu_2}      \int_1^\infty J^{-}_{3\nu_3 }\left(2\pi |y_1|^{1/2}\sqrt{u^2-1}\right) J^{-}_{3\nu_3}\left(2\pi |y_2|^{1/2}\sqrt{1-u^{-2}}\right) u^{3\mu_2} \frac{du}{u},
 \end{equation}
\begin{equation}\label{BesselJ3}
 \mathcal{J}_{3}(y; \mu) = 
 \Bigl| \frac{y_1}{y_2}\Bigr|^{\frac{1}{2}\mu_2}      \int_0^\infty \tilde{K}_{3\nu_3 }\left(2\pi |y_1|^{1/2}\sqrt{1+u^2}\right) J^{-}_{3\nu_3}\left(2\pi |y_2|^{1/2}\sqrt{1+u^{-2}}\right) u^{3\mu_2} \frac{du}{u},
 \end{equation}
\begin{equation}\label{BesselJ4}
 \mathcal{J}_{4}(y; \mu) = 
  \Bigl| \frac{y_1}{y_2}\Bigr|^{\frac{1}{2}\mu_2}      \int_0^1 \tilde{K}_{3\nu_3 }\left(2\pi |y_1|^{1/2}\sqrt{1-u^2}\right) \tilde{K}_{3\nu_3}\left(2\pi |y_2|^{1/2}\sqrt{u^{-2}-1}\right) u^{3\mu_2} \frac{du}{u},
 \end{equation}
   \begin{equation}\label{BesselK}
 \mathcal{J}_{5}(y;\mu) =  \Bigl|\frac{y_1 }{y_2 }\Bigr|^{\frac{1}{2}\mu_2}   \int_0^\infty \tilde{K}_{3\nu_3 }\left(2\pi |y_1|^{1/2}\sqrt{1+u^2}\right) \tilde{K}_{3\nu_3}\left(2\pi |y_2|^{1/2}\sqrt{1+u^{-2}}\right) u^{3\mu_2} \frac{du}{u}.
\end{equation}
Then \cite[Lemma 5]{BB} gives us
\begin{equation}\label{lem++}
   K_{w_6}^{++}(y; \mu)    = \frac{1}{12\pi^2}\frac{ \cos\left(\frac{3}{2}\pi \nu_1\right) \cos\left(\frac{3}{2}\pi \nu_2\right)}{ \cos\left(\frac{3}{2}\pi \nu_3\right)} \mathcal{J}_{5}(y;\mu)
   \end{equation}
   $(y_1, y_2 > 0)$; 
   \begin{equation}\label{lem+-}
   \sum_{w \in \{I, w_4, w_5\}}  K_{w_6}^{+-}(y; w(\mu))    = \frac{1}{24\pi^2} \sum_{w \in \{I, w_4, w_5\}} \Bigl( \mathcal{J}_2(y; w(\mu)) +  \mathcal{J}_3(y; w(\mu)) +  \mathcal{J}_4(y; w(\mu))  \Bigr) ; 
      \end{equation}
 ($y_1 > 0 > y_2$); 
   \begin{equation}\label{lem-+}
  K_{w_6}^{-+}((y_1, y_2); \mu)    =    K_{w_6}^{+-}((y_2, y_1); w_4(-\mu));
         \end{equation}
 ($y_2 > 0 > y_1$); and 
     \begin{equation}\label{lem--}
   \sum_{w \in \{I, w_4, w_5\}} K_{w_6}^{--}(y; w(\mu))  =   \frac{1}{48\pi^2} \sum_{w \in \{I, w_4, w_5\}} \Bigl(4 \mathcal{J}_1^{-}(y; w(\mu)) +2  \mathcal{J}_1^+(y; w(\mu))  \Bigr)
 \end{equation}
($y_1, y_2 < 0$). \\ 

Formulas for $K_{w_4}$ are given in \cite[Definition 1 or Lemma 4]{BB}.  For our purpose it suffices to know that the proof of \cite[Lemma 8]{BB} (which holds a fortiori with the present choice of $h$ in \eqref{defh}) shows that  $\Phi_{w_4}(y)$ and $\Phi_{w_5}(y)$   are negligible (i.e.\ $\ll_B T^{-B}$ for any $B > 0$) if $|y|  \leq T^{3-\varepsilon}$. This is the analogue of \eqref{unless} for the $w_4$, $w_5$ terms, which is already so strong that no further analysis will be necessary.

\section{Analytic preparation} 
In this section we compile a bit of classical analysis that we need in the following. 
For fixed $\sigma \in \Bbb{R}$, real $|t| \geq 10$ and any $M > 0$ we have Stirling's formula
\begin{equation}\label{stir}
  \Gamma(\sigma + it) = e^{-\frac{\pi}{2}|t|} |t|^{\sigma-\frac{1}{2}} \exp\left(i t \log \frac{|t|}{e}\right)g_{\sigma, M}(t) + O_{\sigma, M}(|t|^{-M}),
\end{equation}
where
$$|t|^j \frac{d^j}{d t^j} g_{\sigma, M}(t) \ll_{j, \sigma, M} 1.$$
for all fixed $j \in \Bbb{N}_0$. 

 We record the formula  (see \cite[Section 4.4]{BB})
\begin{equation}\label{intrep1}
\begin{split} 
 \tilde{K}_{it}(x) =   \int_{-\infty}^{\infty} \cos(2x \sinh v)\exp( it v) dv 
\end{split}
 \end{equation}
 for $t \in \Bbb{R}$, $x > 0$. 
 This   integral is not absolutely convergent, but integration by parts shows that we can restrict   (smoothly) to $v  = \pm \log |t|/x + O(1)$ up to a negligible error (see \cite[Section 4.4]{BB} for the precise argument).  From \cite[(4.13)]{BB} we quote
 \begin{equation}\label{besselbound}
  \frac{\partial^j}{\partial x^j} \tilde{K}_{it}(x), \quad \frac{\partial^j}{\partial x^j} J^{\pm}_{ it}(x) \ll_j \left( 1 + \frac{|t|}{x} \right)^{j}
\end{equation}
for $t \in \Bbb{R}$, $x > 0$, $j \in \Bbb{N}_0$.  Moreover, we have the uniform asymptotic expansions 
\begin{equation}\label{unifasymp}
\tilde{K}_{it}(x/2) =  \Re \left( e^{i   \omega(x, t)}  f_M(x, t)\right)  
+ O(|t|^{-M}), \quad \omega(x, t) =  |t| \cdot \text{arccosh} \frac{|t|}{x} - \sqrt{t^2 - x^2},
\end{equation}
for $t \in \Bbb{R}$, $|t| > 1$, $\frac{1}{10} |t|\geq  x > 0$ and fixed $M > 0$ with
\begin{equation}\label{bessdiff1}
x^i |t|^j \frac{\partial^i}{\partial x^i}  \frac{\partial^j}{\partial t^j} f_M(x, t) \ll_{i, j, M}  |t|^{-1/2}
\end{equation}
for any $i, j \in \Bbb{N}_0$; and 
analogously 
\begin{equation}\label{unifasymp1}
 J^{\pm}_{it}(x/2) = 
 \Re \left( e^{i  \tilde{\omega}(x, t)} \tilde{f}_M^{\pm}(x, t)\right) 
 + O(|t|^{-M}),  
\quad \tilde{\omega}(x, t) =  |t| \cdot \text{arcsinh} \frac{|t|}{x} - \sqrt{t^2 + x^2},
\end{equation}
for $t \in \Bbb{R}$,  $|t| > 1$, $x > 0$ and fixed $M > 0$ with
\begin{equation}\label{bessdiff2}
x^i |t|^j \frac{\partial^i}{\partial x^i}  \frac{\partial^j}{\partial t^j} \tilde{f}^{\pm}_M(x, t) \ll_{i, j, M}  \frac{1}{x^{1/2} + |t|^{1/2}}
\end{equation}
for any $i, j \in \Bbb{N}_0$. Again we refer to \cite[Section 4.4]{BB} for details and references\footnote{The formulas after  \cite[(4.17) and (4.20)]{BB}  are stated for $i=0$, but the quoted formulas yield the bounds \eqref{bessdiff1} and \eqref{bessdiff2} for arbitrary $i \in \Bbb{N}_0$.}.

Finally  we quote  two  integration-by-parts lemmas  
from  \cite[Lemma 6]{JM}  and \cite[Corollary 8.3]{BKY}.

\begin{lemma}\label{int-simple} Let $w$ be a smooth function with support in an interval $[a, b]$. Let $r > 0$.  Let $\phi$ be a function that is real-valued on $[a, b]$ and holomorphic on $\mathcal{D} := \{ z \in \Bbb{C} : \inf\{ |z - t| : t \in [a, b]\} < r\}$.  Assume 
\begin{displaymath}
  w^{(j)}(t)  \ll XR^{-j}, \quad   |\phi'(z)| \gg B
\end{displaymath}
for all $t\in [a, b]$, $z \in \mathcal{D}$, $j \in \Bbb{N}_0$ and certain positive real numbers $B, X, R$. Then
\begin{displaymath}
  \int_{-\infty}^{\infty} w(t) e^{  i\phi(t)} dt \ll_j X\left(\frac{1}{B
  \min(R, r)}  \right)^j(b-a) 
\end{displaymath}
for any $j \in \Bbb{N}_0$ with an implied constant depending only on $j$. 
\end{lemma}

\begin{lemma}\label{statphase}
Let $0 < \delta < 1/10$, $X, W, R,   Q > 0$, $Z := Q + X + W + R +1$,  and assume that
\begin{equation*}\label{importantconditions}
Y \geq Z^{3 \delta}, \quad   R \geq \frac{QZ^{ \frac{\delta}{2}} }{W^{1/2}}.\end{equation*} Suppose that $w$ is a smooth function on $\Bbb{R}$ with support on an  interval $J$ of length $R$, satisfying
\begin{equation*}
w^{(j)}(t) \ll_j X R^{-j}
\end{equation*}
for all $j \in \Bbb{N}_0$. Suppose $\phi$ is a real smooth function on $J$ such that there exists a unique point $t_0 \in J$ such that $\phi'(t_0) = 0$, and furthermore
\begin{equation*}
\phi''(t) \gg W Q^{-2}, \quad \phi^{(j)}(t) \ll_j W Q^{-j}, \qquad \text{for } j=2, 3, \dots, t \in J. 
\end{equation*}
Then  there exists a function $w_0(t)$ supported on the interval $[-1, 1]$ such that  \begin{equation*}
\label{eq:shortintegral}
 \int_{\Bbb{R}} w(t) e^{i\phi(t)} dt = \int_{\Bbb{R}} w(t) w_0\Bigl(\frac{t-t_0}{Z^{\varepsilon} |\phi''(t_0)|^{-1/2}}\Bigr) e^{i\phi(t)} dt + O_{A,\varepsilon}(Z^{-A}).
 \end{equation*}
 \end{lemma}

\section{Analysis of the integral transform}\label{analysis}

In this section we analyze the transform $\Phi_{w_6}(y)$ defined in \eqref{defPhi} with a real-analytic test function $\tilde{h}$ satisfying \eqref{diffh} and \eqref{volume} and $K^{\pm, \pm}_{w_6}$ as in \eqref{defK} and \eqref{lem++} -- \eqref{lem--}.  
Ultimately, we will be interested in the double Mellin transform 
\begin{equation}\label{MellinPhi}
\widehat{\Phi}(s, Y) = \sum_{\epsilon_1, \epsilon_2}\Biggl|  \int_{0}^{\infty} \int_{0}^{\infty} g\left(\frac{y_1}{Y_1}, \frac{y_2}{Y_2}\right) \Phi_{w_6}(\epsilon_1 y_1, \epsilon_2 y_2) y_1^{s_1}y_2^{s_2} \frac{dy_1\, dy_2}{y_1y_2} \Biggr|,
\end{equation}
where $Y = (Y_1, Y_2)$ is a pair of   parameters satisfying 
$\log Y_1, \log Y_2 \ll \log T$, 
 $g $ is a  smooth fixed weight function with support in $[1, 2]^2$, say, and $\Re s_1 = \Re s_2 = 0$. We will refer to an error term or a function as being \emph{negligible} if it is $O_B (T^{-B})$ for all $B \geq 0$. It will suffice to study $\widehat{\Phi}(s, Y)$ only for parameters $Y_1, Y_2$ where $\Phi_{w_6}(y_1, y_2) g(y_1/Y_1, y_2/Y_2)$ is not negligible, in particular we will always assume \eqref{unless}.

By Mellin inversion and \eqref{defK} we have
\begin{equation}\label{doublemellin}
\begin{split}
\widehat{\Phi}(s, Y) = \sum_{\epsilon_1, \epsilon_2}\Bigl|  \int_{\Re \mu = 0} & \tilde{h}(\mu) \int_{\Re u = \varepsilon} (4\pi^2)^{-u_1-u_2} G(u, \mu) S^{\epsilon_1, \epsilon_2}(u; \mu)\\
& \times \widehat{g}(s-u)  Y_1^{s_1-u_1} Y_2^{s_2-u_2} \frac{du}{(2\pi i)^2} \text{spec}(\mu) d\mu \Bigr|.
\end{split}
\end{equation}
We write $$\Im u_j = v_j, \quad \Im s_j = t_j, \quad   \Im \mu_j = \tau_j$$
and recall that $\Re s_j = 0$ (and the integration in \eqref{doublemellin} is over $\Re u_j = \varepsilon$, $\Re \mu_j = 0$). We keep these lines fixed and will not shift any contours.  By the rapid decay of $\widehat{g}$ on vertical lines this multiple integral is absolutely convergent, and  outside the range  $v_j =  t_j + O(T^{\varepsilon})$ it is negligible. 

We can also express $\widehat{\Phi}(s, Y)$ in terms of  the integral representations \eqref{BesselJ1} -- \eqref{BesselK}, and we make some preliminary comments.  In view of \eqref{lem++}--\eqref{lem--} and the Weyl symmetry of $\tilde{h}$, it suffices to analyze 
\begin{equation}\label{mellin-j}
  \widehat{\Phi}_j(s, Y) :=  \int_{0}^{\infty} \int_{0}^{\infty}  \Phi_{j}(y_1, y_2) y_1^{s_1}y_2^{s_2} \frac{dy_1\, dy_2}{y_1y_2}
\end{equation}
for $1 \leq j \leq 5$, where 
\begin{equation}\label{Phij}
\Phi_j(y, Y) := g\left(\frac{y_1}{Y_1}, \frac{y_2}{Y_2}\right) \int_{\Re \mu = 0}  \tilde{h}(\mu) \mathcal{J}_j(y; \mu) \text{spec}(\mu)d\mu. 
\end{equation}
(In particular, as in the proof of \cite[Lemma 9]{BB}, in the case of $\Phi_5$ we can restrict ourselves to the positive Weyl chamber, where $\cos(\frac{3}{2} \pi\nu_1)\cos(\frac{3}{2} \pi\nu_2) \cos(\frac{3}{2} \pi\nu_3)^{-1} = 1/2 + O(T^{-B})$, so this factor can be disregarded at the cost of a negligible error.)  

 In the $\mu$-plane we will always use the coordinates 
 \begin{equation}\label{coordinates}
 \tau_2 = \Im \mu_2  \quad \text{and} \quad \rho = 3\Im\nu_3,
 \end{equation}
 and   the support of $\tilde{h}$ restricts $\rho$ to an interval $I$ (say) of length $\ll T^{1-\varepsilon}$, where $|\rho| \asymp T$. By symmetry we will assume with loss of generality that $\rho > 0$. 
 Notice that 
 $$\frac{1}{2}(\rho - \tau_2) = \Im \mu_1, \quad -\frac{1}{2}(\rho + \tau_2) = \Im \mu_2.$$
 When we integrate by parts using Lemma \ref{int-simple}  with respect to $\rho$, we need to bound derivatives of various phase functions in a complex disc about $I$ that we denote by 
 $$D_{\rho} := \{ z \in \Bbb{C} : \inf\{|z - t| : t \in I\} < T^{1-\varepsilon}\}.$$\\

Our first aim is to show that we can effectively truncate $t_1$ and $t_2$, so that from now on  all parameters under consideration are at most powers of $T$. 
 The exponent $T^5$ in the following bound comes from \eqref{volume}. It can easily be improved (and we will do so in the following lemmas); the point of this lemma is only the truncation of $t_1, t_2$.

\begin{lemma}\label{phi1} We have $$\widehat{\Phi}(s, Y) \ll_B T^{5}\left(\frac{|t_1|}{T+ Y_1^{1/2} + Y_1^{1/3} Y_2^{1/6}} + \frac{|t_2|} {T+ Y_2^{1/2} + Y_2^{1/3} Y_1^{1/6}}\right)^{-B}$$ 
for all $B \geq 0$. 
\end{lemma}

\textbf{Proof.} We show the desired bound for each $\widehat{\Phi}_j(s, Y)$. It follows from \eqref{besselbound} and repeated integration by parts that in general 
$$\int_{-\infty}^{\infty}  y^{\pm \frac{1}{2} i \tau_2} J^{\pm}_{3\nu_3} (y^{1/2} \alpha) g\left(\frac{y}{Y}\right) y^{it_1}  \frac{dy}{y} \ll_B \left(\frac{Y^{1/2} \alpha + T}{|t_1|}\right)^B$$
for $\tau_2, \rho \ll T$ and $\alpha > 0$, and the same holds for  $\tilde{K}_{3\nu_3}$ in place of $J^{\pm}_{3\nu_3}$. 
Integrating \eqref{BesselJ1} -- \eqref{BesselK} over   $\tau_2$ restricts $u$ to $u \asymp (y_2/y_1)^{1/6}$, up to a negligible error. Applying the previous bound with $\alpha = \sqrt{|1 \pm u^{\pm 2}|} \ll 1 + u^{\pm 1}$,  we can save arbitrary powers of $|t_1|$ (respectively $|t_2|$)  once  $|t_1|$ grows beyond $ T+Y_1^{1/2} + Y_1^{1/3} Y_2^{1/6} $ (respectively  $|t_2|$ grows beyond $ T+Y_2^{1/2} + Y_2^{1/3} Y_1^{1/6}$).  
\hfill $\square$\\

\begin{lemma}\label{phi2}
Let $T_1, T_2 \geq T^{1+\varepsilon}$ be two parameters, and let
$$S_1 := \sup_{|t_1| \asymp T_1} \int_{|t_2| \asymp T_2}  |\widehat{\Phi}((it_1, it_2), Y)| dt_2, \quad S_2 := \sup_{|t_2| \asymp T_2}  \int_{|t_1| \asymp T_1}  |\widehat{\Phi}((it_1, it_2), Y)| dt_1.$$
Then 
\begin{equation}\label{s1s2a}
(S_1S_2)^{1/2} \ll T^{2+\varepsilon}.
\end{equation}
and
\begin{equation}\label{s1s2b}
(S_1S_2)^{1/2} \ll  \frac{T^{5+\varepsilon}(T_1+T_2)^{1/2}}{T_1T_2}.
\end{equation}
 \end{lemma}

\textbf{Proof.} For notational simplicity let us assume (by symmetry) that $T_1 \geq T_2$. We start by bounding $\widehat{\Phi}((it_1, it_2), Y)$ under the assumption $|t_1| \asymp T_1$, $|t_2| \asymp T_2$ together with the additional assumption $|t_1+ t_2| \geq T^{\varepsilon}$. 
By \eqref{stir}, the oscillation of the $\mu$-integral in \eqref{doublemellin} has the phase 
$$\phi(\mu) := \sum_{j=1}^3 (v_1 -\tau_j) \log \frac{|v_1 - \tau_j|}{e} +  (v_2 + \tau_j) \log \frac{|v_2 + \tau_j|}{e}.$$ 
If   $v_j = t_j + O(T^{\varepsilon})$, 
 our assumption $\min(T_1, T_2)  \geq T^{1+\varepsilon}$ implies that $|v_1|, |v_2| \gg T^{1+\varepsilon}$ dominate all $\tau_j \ll T$, and we conclude by the mean value theorem  that (recall our choice of coordinates \eqref{coordinates}) 
$$\Bigl|\frac{\partial}{\partial \rho} \phi(\mu)\Bigr| =\left| \frac{1}{2} \log \Bigl| \frac{(v_1 +  \frac{1}{2}(\tau_2 + \rho))(v_2 - \frac{1}{2}(\tau_2 - \rho))}{(v_1 + \frac{1}{2}(\tau_2 - \rho))(v_2 - \frac{1}{2}(\tau_2 + \rho))}\Bigr|\right| \asymp \Bigl|\rho\left(\frac{1}{v_1} + \frac{1}{v_2}\right)\Bigr| \asymp \frac{T |v_1+v_2|}{T_1T_2} \asymp \frac{T |t_1+t_2|}{T_1T_2}. $$
 By \eqref{diffh} and Lemma \ref{int-simple}  with $R = r = T^{1-\varepsilon}$, we conclude that the $\mu$-integral is negligible, unless
\begin{equation}\label{mu-int}
 (1 +  |t_1 + t_2|) T^{2-2\varepsilon} \ll T_1T_2.
 \end{equation}
Since $T_1 \geq T_2 \geq T^{1+\varepsilon}$, this remains trivially true if  $|t_1+ t_2| \leq T^{\varepsilon}$, so that \eqref{mu-int} holds in all cases. 

Now we estimate \eqref{doublemellin}  trivially using \eqref{volume} and Stirling's formula for the function $G(u, \mu)$ defined in \eqref{Gsmu} and obtain 
\begin{equation}\label{trivial}
\widehat{\Phi}(s, Y) \ll T^{\varepsilon}  \frac{T^5(1+ |t_1 + t_2|)^{1/2}}{T_1^{3/2} T_2^{3/2}}.   
\end{equation}
This gives immediately \eqref{s1s2b}. 

 If $T_1 \geq C T_2$ for a sufficiently large constant $C$ so that $|t_1+t_2|  \asymp T_1$, then \eqref{mu-int} implies $T_1 \geq T_2 \gg T^{2-2\varepsilon}$, and together with \eqref{trivial} we conclude  \eqref{s1s2a}. On the other hand, if $T_1 \asymp T_2$, then \eqref{mu-int} and the previous bound imply that, up to a negligible error, we have
\begin{displaymath}
\begin{split}
S_1 &\ll T^{\varepsilon} \sup_{|t_1| \asymp T_1} \int_{t_2 = -t_1 + O(T_1T_2 T^{2\varepsilon - 2})} \frac{T^5(1+ |t_1 + t_2|)^{1/2}}{T_1^{3/2} T_2^{3/2}} \ll    \frac{T^{5+\varepsilon}  (T_1T_2/T^2)^{3/2}}{(T_1T_2)^{3/2}}
\end{split}
\end{displaymath}
and similarly for $S_2$, and \eqref{s1s2a} follows again. \hfill $\square$\\

The analysis in the previous proof was greatly simplified by the assumption $T_1, T_2 \geq T^{1+\varepsilon}$, so that the $v_j$ and the $\tau_j$ could not interfere. Our final lemma in this section complements the previous two lemmas and shows in particular  that \eqref{s1s2a} remains true without the extra assumption $T_1, T_2 \geq T^{1+\varepsilon}$. 

\begin{lemma}\label{phi3} a) Suppose that $\min(T_1, T_2) \leq T^{1+\varepsilon}$. Then $\widehat{\Phi}(s, Y)$ is negligible for $|t_1| \asymp T_1$ and $|t_2| \asymp T_2$ unless $Y_1 \asymp Y_2$, in which case 
  \eqref{s1s2a} holds.\\ 
b) Suppose that $Y_1 \geq CY_2$ or $Y_2 \geq CY_1$ for a sufficiently large constant $C$. Then 
$\widehat{\Phi}(s, Y)$ is negligible unless
\begin{equation}\label{t-lower}
|t_1| \gg (Y_1^{1/2} + Y_1^{1/3} Y_2^{1/6})T^{-\varepsilon} \gg T^{2-\varepsilon}, \quad |t_2| \gg  (Y_2^{1/2} + Y_2^{1/3} Y_1^{1/6})T^{-\varepsilon} \gg T^{2-\varepsilon}.
\end{equation}
\end{lemma}

\textbf{Remark.} The assumptions in a) and b) look a bit artificial, and indeed we could make much more general statements at the cost of more work. The statement of this lemma is tailored precisely to our needs in the proof of Theorem \ref{thm2}. \\

\textbf{Proof.} We analyze $\Phi_j(y, Y)$  
and $\widehat{\Phi}_j(s, Y)$, defined in \eqref{mellin-j} and \eqref{Phij}, for $j = 1, \ldots, 5$. 
We will always denote the argument of the two Bessel functions in \eqref{BesselJ1} -- \eqref{BesselK} by $x_1, x_2$:
$$x_1 = 4\pi y_1^{1/2} \sqrt{|1 \pm u^2|}, \quad   x_2 = 4\pi y_2^{1/2} \sqrt{|1 \pm u^{-2}|}$$
with appropriate sign depending on $j \in \{1, \ldots, 5\}$. By \eqref{diffh}, the $\tau_2$-integral   restricts us (up to a negligible error) to 
 \begin{equation}\label{usize}
u = \left(\frac{y_2}{y_1}\right)^{1/6} (1 + O(T^{\varepsilon-1})) \asymp \left( \frac{y_2}{y_1}\right)^{1/6}.
\end{equation}
We always think of this region as extracted smoothly. Before we turn to the individual cases, we make some general comments. Suppose we can show that some portion $\Psi(y, Y)$ of some $\Phi_j(y, Y)$  is negligible unless 
\begin{equation}\label{cond1}
   y_1 = y_2 \left(1 + O (T^{\varepsilon - 1})\right)
    \end{equation} 
(so that in particular $Y_1 \asymp Y_2 =: Y_0$, say),   in which case
\begin{equation}\label{diff1}
  y^j  \frac{\partial^j}{\partial y^j}  \Psi\bigl((y,  y + z), Y\bigr) \ll_{j } T^{3  + \varepsilon}
  \end{equation}
  for all $j \in \Bbb{N}_0$ and $z \ll Y_0T^{\varepsilon - 1}$. Then by trivial estimates and sufficiently many integrations by parts we have
$$\widehat{\Psi}(s, Y) \ll T^{2+\varepsilon} \left(1 + \frac{|s_1 + s_2|}{T^{\varepsilon}}\right)^{-B}$$
for $\Re s_1 = \Re s _2 = 0$, and so 
$$\Biggl(\sup_{|t_1| \asymp T_1} \int_{|t_2| \asymp T_2}  |\widehat{\Psi}(s, Y )| dt_2 \cdot  \sup_{|t_2| \asymp T_2}  \int_{|t_1| \asymp T_1}  |\widehat{\Psi}(s, Y )| dt_1\Biggr)^{1/2} \ll T^{2+\varepsilon},$$
which is a version of \eqref{s1s2a}. 

In each of the cases $j=1, \ldots 5$ we will split $\Phi_j(y, Y)$ into various pieces depending on signs and sometimes size conditions of various parameters, and we will show that each piece is 
\begin{enumerate}
\item\label{1} negligibly small; or 
\item\label{2}   \eqref{cond1} and \eqref{diff1} hold, in which case part a) of the lemma follows and part b) is void; or 
\item\label{3}   the double Mellin transformation is negligible unless \eqref{t-lower} holds, in which case part b) of the lemma holds and part a) is void; or
\item\label{4}   $Y_1 \asymp Y_2$ holds and its double Mellin transform is negligible unless $|t_1|, \, |t_2| >T^{1+\varepsilon}$, in which case we have nothing to prove, because both parts of the lemma are void. 
\end{enumerate}
To this end we insert either the uniform asymptotic expansion \eqref{unifasymp1} or the integral representation \eqref{intrep1} and apply integration by parts in the form of Lemma \ref{int-simple} (or sometimes Lemma \ref{statphase}) in the $\rho$ or $u$-integral in \eqref{Phij}  or in $y_1, y_2$-integral in \eqref{mellin-j}.  Note that we have already squeezed out all information from the $\tau_2$-integral in \eqref{usize}. \\


\textbf{The case $j=1$.} We insert the uniform asymptotic formula \eqref{unifasymp1}. This expresses $\Phi_1(y, Y)$, up to a negligible error, as a sum of 
\begin{displaymath}
\begin{split}
 \Psi_{\epsilon_1, \epsilon_2}(y, Y) =  g\left(\frac{y_1}{Y_1}, \frac{y_2}{Y_2}\right) \int_{(0)}  \tilde{h}(\mu)   
  \int_0^\infty & e^{i\phi(u, \rho; y_1, y_2)}  \tilde{f}_M^{\pm}(x_2, \rho)\tilde{f}_M^{\pm}(x_1, \rho)
   \frac{du}{u} \, \text{spec}(\mu)  d\mu, 
 \end{split}
\end{displaymath}
where
$$\phi(u, \rho; y_1, y_2) =  \epsilon_1  \tilde{\omega} (x_2, \rho )  + \epsilon_2   \tilde{\omega} (x_1, \rho ) + 3\tau_2 \log u + \frac{1}{2} \tau_2 \log \frac{y_1}{y_2}$$
and 
$\epsilon_1, \epsilon_2 = \pm 1$ (which are independent of the $\pm$ signs of the weight functions $f_M^{\pm}$). Without loss of generality we assume $Y_1 \gg Y_2$ (by symmetry), in particular $u \ll 1$.  
Then by \eqref{usize} we have 
\begin{equation}\label{size-x1x2}
x_1 \asymp Y_1^{1/2} \gg x_2 \asymp  Y_2^{1/3} Y_1^{1/6}.
\end{equation}
We compute
\begin{equation}\label{u}
  \frac{\partial}{\partial u }\phi(u, \rho; y_1, y_2) =    \epsilon_1  \frac{\sqrt{\rho^2  + x_2^2}}{u(1+u^2)}- \epsilon_2  \frac{u \sqrt{\rho^2  + x_1^2}}{ 1+u^2}  + \frac{3\tau_2}{u}
  \end{equation} 
and 
\begin{equation}\label{rho}
  \frac{\partial}{\partial \rho }\phi(u, \rho; y_1, y_2) = \epsilon_1 \text{arcsinh}\left(\frac{\rho}{x_2}\right) + \epsilon_2 \text{arcsinh}\left(\frac{\rho}{x_1}\right). 
\end{equation}
 
We now distinguish two cases. If $\epsilon_1 = \epsilon_2 = \epsilon$ (say), then applying Lemma \ref{int-simple} with $R = r = T^{1-\varepsilon}$, $B = T^{2\varepsilon - 1}$, and noting that\footnote{Indeed, we have $\arg  [\partial_b  \text{arcsinh}(a+ib) ] =  \arg[ i (1 + (a+ib)^2)^{-1/2}] \in [0, \pi/2]$ for $a, b > 0$, so that  $\Re\, \text{arcsinh}(a+ib)$ is increasing in $b$ for $a, b > 0$, and analogously  decreasing for $a > 0 > b$.} $\Re\, \text{arcsinh}(a+ib) \geq \text{arcsinh}(a)$ for $a > 0$, $ b \in \Bbb{R}$, we see by \eqref{diffh}, \eqref{bessdiff2} and  \eqref{rho}  that the $\rho$-integral is negligible unless both
$\rho/x_2, \rho/x_1 \ll T^{2\varepsilon - 1}$, or in other words, 
\begin{equation}\label{big}
  Y_2^{1/3} Y_1^{1/6}  \gg T^{2-2\varepsilon}. 
\end{equation}
With this information we consider now the double Mellin transform  $\widehat{\Psi}_{\epsilon, \epsilon}(s, Y)$ and show that $t_1$, $t_2$ must be large. Indeed, we have
\begin{equation}\label{indeed}
\Bigl|\frac{\partial}{\partial y_1}\Bigl( \phi(u, \rho; y_1, y_2)  + t_1 \log y_1 + t_2 \log y_2\Bigr) \Bigr| = \Bigl|\epsilon \frac{\sqrt{\rho^2 + x_1^2}}{2y_1} + \frac{t_1}{y_1} + \frac{ \tau_2}{2y_1}\Bigr| \gg \frac{x_1}{Y_1} \gg \frac{1}{Y_1^{1/2}}
\end{equation}
(in which case the $y_1$-integral is negligible by Lemma \ref{int-simple}) unless $|t_1| \asymp  x_1 \asymp Y_1^{1/2} \gg T^{2-\varepsilon}$ by \eqref{big}, and a similar relation holds for $t_2$. We are therefore in situation \eqref{3}.

If $\epsilon_1 \not= \epsilon_2$, we consider the derivative with respect to $u$. If $Y_1^{1/6} Y_2^{1/3} \gg T^{1+\varepsilon},$ then   \eqref{size-x1x2} and \eqref{usize} imply that \eqref{u} is $\asymp Y_1^{1/3}Y_2^{1/6} \gg T^{1+\varepsilon}$ in a disc given by \eqref{usize}, hence applying Lemma \ref{int-simple} with $B = T^{1+\varepsilon}$, $R =r  = T^{\varepsilon - 1}$,  
 we see   that the $u$-integral is negligible. Thus we conclude 
\begin{equation}\label{conclude}
Y_1^{1/6} Y_2^{1/3} \ll T^{1+\varepsilon}.
\end{equation}
However, by Lemma \ref{int-simple} the $\rho$-integral still forces  
$\partial_{\rho} \phi(u, \rho; y_1, y_2) \ll T^{2\varepsilon - 1}$ somewhere in $D_{\rho}$. Since \eqref{unless} and \eqref{conclude} imply $Y_1 \asymp Y_2 =: Y_0 = T^{2+o(1)}$ and hence $x_1, x_2 = T^{1+o(1)}$, the mean value theorem and \eqref{usize}  imply 
\begin{displaymath}
\begin{split}
T^{2\varepsilon - 1}& \gg \frac{|x_1 - x_2| T}{x_1(T + x_1)} \gg\frac{ |x_1 - x_2| }{T^{1+\varepsilon}} = 4\pi\sqrt{1 + u^2}   \frac{| y_1^{1/2} - u^{-1} y_2^{1/2}|}{T^{1+\varepsilon}}  \asymp  \frac{|y_1^{1/2} - u^{-1} y_2^{1/2}|}{T^{1+\varepsilon}} \\
&= \frac{|y_1^{1/2} -   y_2^{1/2}| + O(T^{\varepsilon})}{T^{1+\varepsilon}} . 
\end{split}
\end{displaymath}
We conclude  (re-defining $\varepsilon$)
\begin{equation}\label{stronger}
T^{2-\varepsilon} \ll y_1 = y_2 (1 + O(T^{\varepsilon-1})).
\end{equation}

For $0 < \alpha \ll 1$ and   $y \in \Bbb{C}$ with $\Re y \asymp Y_0$, $\Im y \ll Y_0$   we have by direct computation and two applications of the mean value theorem that
\begin{equation}\label{partial-y}
\begin{split}
\frac{\partial}{\partial y}  \big( \tilde{\omega}(y^{1/2} \alpha, \rho) - \tilde{\omega}((y+z)^{1/2} \alpha u^{-1}, \rho)\big) &= \frac{1}{2} \left(\frac{\sqrt{\rho^2 + \alpha^2 (y+z)/u}}{y+z} - \frac{\sqrt{\rho^2 + \alpha^2 y}}{y} \right)\\
& \ll \frac{(\rho + \alpha Y_0^{1/2})|z|}{Y_0^2} + \frac{\alpha^2 |u-1|}{\rho+ \alpha Y_0^{1/2}},  
\end{split}
\end{equation}
provided that 
$|z| \leq cY_0$ for a sufficiently small constant $c > 0$ and 
$u \asymp 1$. 
By Cauchy's integral formula for higher derivatives on a circle of radius $\asymp Y_0$ we conclude with $\alpha = 4\pi \sqrt{1 + u^2} \asymp 1$ for $y \asymp Y_0$ that 
\begin{equation}\label{osc}
\begin{split}
     y^j  \frac{\partial^j}{\partial y^j}    \phi (u, \rho; y, y+z )   \ll_j  \frac{|z| \sqrt{\rho^2 + Y_0}}{Y_0} + \frac{|u-1| Y_0}{\sqrt{\rho^2 + Y_0}}  + \frac{|\tau_2| |z|}{Y_0} \ll T^{\varepsilon}
     \end{split}
     \end{equation}
for all $j \in \Bbb{N}_0$ and  $$y \asymp Y_0, \quad z \ll Y_0 T^{-1+\varepsilon}, \quad T^{2-\varepsilon} \ll Y_0  \ll T^{2+\varepsilon}, \quad  |u-1| \ll  \frac{|z|}{y} + \frac{1}{T^{1-\varepsilon} } \ll  \frac{1}{T^{1-\varepsilon} }$$
(the first inequality for $u$  holds   by \eqref{usize}). Hence also 
\begin{equation*} 
\begin{split}
     y^j  \frac{\partial^j}{\partial y^j}   e^{i  \phi (u, \rho; y, y+z )}   \ll_j  \frac{|z| \sqrt{\rho^2 + Y_0}}{Y_0} + \frac{|u-1| Y_0}{\sqrt{\rho^2 + Y_0}}  + \frac{|\tau_2| |z|}{Y_0} \ll T^{\varepsilon}
     \end{split}
     \end{equation*}
under the same conditions. 
By \eqref{volume}, \eqref{usize} (which restricts $u$ to an interval of length $T^{\varepsilon-1})$  and \eqref{bessdiff2}  
we conclude that
$$y^j  \frac{\partial^j}{\partial y^j} \Psi_{\epsilon, -\epsilon}((y, y+z), Y) \ll T^{3+\varepsilon} $$
so that we are in situation \eqref{2}. \\ 

\textbf{The case $j= 5$.}   
  We may continue to assume $Y_1 \gg Y_2$ (so that $u \ll 1$).  
  Moreover, \eqref{unless} together with the rapid decay of the Bessel $K$-function implies 
  \begin{equation}\label{moreover}
  T^{2-\varepsilon} \ll Y_1, Y_2  \ll T^{2}, \quad T^{1-\varepsilon} \ll x_1, x_2 \ll T.
  \end{equation}  We insert the integral representation \eqref{intrep1} for both Bessel functions, expressing them as
\begin{equation*}
\int \exp(i \rho (v-w)) \exp\left(\pm  ix_1 \sinh v \pm   i x_2\sinh w\right) dv \, dw, 
\end{equation*}
where both integrals are understood to be cut off smoothly at $\exp |v|  \ll  T/T^{1-\varepsilon} =   T^{\varepsilon} .$  
 Integrating over $\rho$ shows $|v - w| \ll T^{\varepsilon - 1}$ (up to a negligible error). Hence we may write $w = v + \eta$, $\eta  \ll T^{\varepsilon - 1}$, so that the $v$-integral is given by
$$\int f(v)  \exp\left(\pm i ( x_1 \pm  x_2 \cosh \eta)  \sinh v   \right) dv,$$
where $f$ is a   smooth   weight function with support restricted to $\exp |v| \ll T^{\varepsilon} $ and satisfying $f^{(j)} \ll_{j, \varepsilon} T^{j \varepsilon}$, in which we have absorbed the exponential $\exp(\pm ix_2  \sinh \eta \cosh v)$. Changing variables, this equals
$$ \int  \frac{f({\rm arcsinh} (v))}{\sqrt{1 + v^2}}  \exp (\pm i ( x_1 \pm  x_2 \cosh \eta) v   ) dv,$$
which forces $ x_1 \pm  x_2 \cosh \eta \ll T^{\varepsilon}$, and therefore also  
\begin{equation}\label{x1x2}
   x_1 \pm  x_2 \ll T^{\varepsilon}.
   \end{equation}
    This is only possible if the $\pm$ sign is negative, and (using \eqref{usize}) we conclude as in the proof of the case $j=1$  that \eqref{stronger} must hold. To verify that we are situation \eqref{2}, 
   we notice that we save a factor $T^{1-\varepsilon}$ from both the $\eta$ and the $u$-integral giving us an upper bound $\Phi_5(y, Y) \ll T^{3+\varepsilon}$, and by \eqref{moreover} and \eqref{usize} we have  
    \begin{equation}\label{check}
    \begin{split}
     y^j & \frac{\partial^j}{\partial y^j} \Bigl( \exp\Bigl(\pm ( ix_1 \sinh v -  i x_2\sinh (v+\eta))\Bigr)  \left(1 + \frac{z}{y}\right)^{\frac{1}{2}i\tau_2} \Bigr)\\
&   \ll_j \left(  |\sinh v| \Bigl(\frac{|\tau_2| |z|}{y} + \frac{|z|}{y^{1/2}} + (|u-1| + |\eta|)y^{1/2}\Bigr)\right) \ll T^{\varepsilon}
   \end{split}
    \end{equation}
  for $z \ll y T^{\varepsilon - 1}$, $   \sinh v \ll T^{\varepsilon}$ and $\eta \ll T^{\varepsilon - 1}$. \\
    

\textbf{The case $j=3$.}    This case is not symmetric in $y_1$ and $y_2$, and interestingly it turns out that $\Phi_3(y, Y)$ is always negligible.  
The rapid decay of the Bessel $K$-function and \eqref{unless} imply that we may assume
\begin{equation}\label{bessel}
T^{1-\varepsilon} \ll   x_1 \ll T.
 \end{equation}  
    We insert the integral representation \eqref{intrep1} for the $K$-function and the uniform asymptotic expansion \eqref{unifasymp1} for the $J$-function, so that the phase in the present situation is given by
$$\phi(v, u, \rho) = \epsilon_1 \tilde{\omega} ( x_2, \rho )  +  \rho  v + \epsilon_2 x_1 \sinh v + 3  \tau_2 \log u + \frac{1}{2} \tau_2 \log \frac{y_1}{y_2}$$
with $\exp|v| \ll  T^{\varepsilon}$ as in the case $j=5$ and $\epsilon_1, \epsilon_2 = \pm 1$.  By the symmetry of the Bessel-$K$-function we may assume $\rho > 0$. 
As before we see that the $\rho$-integral is negligible unless
  \begin{equation}\label{rho-diff}
  \Bigl|\frac{\partial}{\partial \rho} \phi(v, u, \rho)\Bigr| = \Bigl|v +\epsilon_1 \text{arcsinh}\left(\frac{\rho}{x_2}\right)\Bigr| \ll T^{2\varepsilon - 1}
  \end{equation}
somewhere in $D_{\rho}$. Let us first assume that 
\begin{equation}\label{assumption}
Y_2^{1/2} \gg T^{1+\delta}
\end{equation}
 for some fixed $0 < \delta < 1/10$ (so that $u \gg 1$ by \eqref{usize}). Then \eqref{rho-diff} implies \begin{equation}\label{vsize} v \ll T^{-\delta}, 
\end{equation}  
and we extract this range smoothly.   Now looking at the $v$-integral, we have
  \begin{equation*}\label{v-diff}
  \frac{\partial}{\partial v} \phi(v, u, \rho) = \rho +\epsilon_2 x_1 \cosh v.
  \end{equation*}
  By Lemma \ref{int-simple}, this must be $\ll T^{\delta+\varepsilon}$  in the disc \eqref{vsize}, otherwise the $v$-integral is negligible.  In particular,  we must have  $\epsilon_2 = -1$, and together with \eqref{usize} and \eqref{bessel} we conclude that we can restrict to 
\begin{equation}\label{rhosize}
\begin{split} 
  \rho  = 4\pi y_2^{1/6} y_1^{1/3}(1 + O(T^{-2\delta/3 + \varepsilon})), 
  \end{split}
  \end{equation}
  where we used \eqref{bessel} several times. Again we extract this range smoothly. With this information we consider the $u$-integral. The derivative of its phase is given by
  \begin{displaymath}
  \begin{split}
 \frac{\partial}{\partial u} \phi(v, u, \rho)  & = \epsilon_1  \frac{\sqrt{\rho^2 u^2 + 16\pi^2y_2 (1 + u^2)}}{u^2(1+u^2)} -\epsilon_2 \frac{4\pi y_1^{1/2} u \sinh v}{\sqrt{1 + u^2}} + \frac{3 \tau_2}{u}\\
  & =  \left(\epsilon_1 \frac{ 4\pi y_2^{1/2}}{u^3} + \frac{3 \tau_2}{u} \right)(1 + O(T^{-\varepsilon})) = \frac{\epsilon_1 \rho + 3 \tau_2}{u} (1 + O(T^{-\varepsilon})) \\
  & = \frac{1}{u}  \left\{ \begin{array}{l} 6\Im \nu_2\\ - 6\Im \nu_1\end{array}\right\}(1 + O(T^{-\varepsilon}))  \asymp \frac{T}{u}
  \end{split}
  \end{displaymath}
 for some sufficiently small $\varepsilon$ in the disc described by \eqref{usize}, where we used  \eqref{bessel}, \eqref{vsize}, \eqref{rhosize} and the fact that $\rho = 3 \Im \nu_3$. 
Hence by Lemma \ref{int-simple} the $u$-integral is negligble under 
the present assumption \eqref{assumption}.   

Together with \eqref{bessel} and \eqref{unless} we may therefore assume $T^{2-\varepsilon} \ll Y_1, Y_2 \ll T^{2+\varepsilon},$ and hence 
$T^{1-\varepsilon} \ll x_1, x_2 \ll T^{1+\varepsilon}.$  
Let us assume $v>0$, the other case is essentially identical.  Then $\epsilon_1 = -1$ by \eqref{rho-diff}. We consider again in more detail the $v$- and the $\rho$-integral, noting that
$$
  \frac{\partial^2}{\partial \rho^2} \phi(v, u, \rho) \asymp \frac{1}{x_2 + T}, \quad   \frac{\partial^j}{\partial \rho^j} \phi(v, u, \rho) \ll \frac{1}{(x_2 + T)^{j-1}}\quad (j\geq 2).$$
Applying Lemma \ref{statphase} with 
$W = Q = x_2 + T$, $R = T^{\varepsilon - 1}$, we can restrict the range of $\rho$ to 
\begin{equation}\label{eq-v}
v = \text{arcsinh}\left(\frac{\rho}{ x_2}\right) + O(T^{\varepsilon - 1/2}). 
\end{equation}
For the moment we only use that this implies $\sinh v \gg T^{-\varepsilon} \cosh v$ and apply the same argument to the $v$-integral. Here we have 
$$\frac{\partial^2}{\partial v^2} \phi(v, u, \rho) =  x_1 \sinh v \gg T^{1-\varepsilon} \quad \frac{\partial^j}{\partial v^j} \phi(v, u, \rho) \ll   x_1 \cosh v \ll T^{1+\varepsilon} \quad (j\geq 2).$$
Applying Lemma \ref{statphase} with 
$W =T^{1-3\varepsilon}$, $Q = T^{-\varepsilon}$, $R = T^{-\varepsilon}$, we can restrict the $v$-integral to 
\begin{equation}\label{eq-rho}
\rho =   x_1 \cosh v + O(T^{1/2 + \varepsilon})
\end{equation}
and we must have $\epsilon_2 = -1$. 
Solving \eqref{eq-v} and \eqref{eq-rho}  for $\rho$ under the present size conditions gives
\begin{equation}\label{solve}
\rho = \frac{x_1x_2}{\sqrt{x_2^2 - x_1^2}} + O(T^{1/2 + \varepsilon}). 
\end{equation}
(In particular, $x_2$ must be a bit larger than $x_1$, otherwise the integral is negligible.)  Once again we consider now the $u$-derivative
\begin{displaymath}
  \begin{split}
 \frac{\partial}{\partial u} \phi(v, u, \rho)  & = -  \frac{\sqrt{\rho^2 u^2 + 16\pi^2y_2 (1 + u^2)}}{u^2(1+u^2)} - \frac{4\pi y_1^{1/2} u \sinh v}{\sqrt{1 + u^2}} + \frac{3 \tau_2}{u},  
 \end{split}
  \end{displaymath}
and we substitute \eqref{usize}, \eqref{eq-v} and \eqref{solve}. After a marathon of simplification (noting that \eqref{solve} implies $|x_1 - x_2| = T^{1/2 + o(1)}$), we arrive at
$$ \frac{\partial}{\partial u} \phi(v, u, \rho) = \frac{-\rho + 3 \tau_2}{u}(1 + O(T^{-1/2+\varepsilon})) \asymp \frac{T}{u},$$
 and as before we use Lemma \ref{int-simple} to show that the $u$-integral is negligible.    \\
 
 \textbf{The case $j=4$.}   Since $u\leq 1$ in this case,   we have necessarily $Y_1 \gg Y_2$. As in the case $j=5$ we insert the integral representation \eqref{intrep1} and express the product of the two Bessel functions as 
 \begin{equation*}
\int \int  \exp(i \rho (v-w)) \exp\left(\pm   i x_1 \sinh v \pm i x_2 \sinh w\right) dv\, dw, 
\end{equation*}
where    the  integrals are restricted (smoothly) to $\exp|v| \ll   T/x_1$, $\exp|w| \ll  \log T/x_2.$ Note that the rapid decay of the Bessel $K$-function restricts to 
\begin{equation}\label{bessel-x}
  x_1, x_2 \ll  T
  \end{equation}
(but  we have a priori no lower bounds).  
 Integrating over $\rho$ shows $|v - w| \ll T^{\varepsilon - 1}$ (up to a negligible error). As in \eqref{x1x2} we conclude that
 \begin{equation}\label{small}
x_1  \pm x_2 \ll T^{\varepsilon}.
 \end{equation}
By \eqref{unless} and \eqref{usize}, this is clearly impossible if $Y_1 \geq C Y_2$ for some sufficiently large constant $C$. 

Let us therefore now assume  $Y_1 \asymp Y_2 =: Y_0 \gg T^{2-\varepsilon}$ (by \eqref{unless}), so that in particular $u \asymp 1$.  Then we make a smooth dyadic decomposition of $u$ close to 1 and   put 
$$1 \gg V \asymp 1 - u, \quad \Delta := |y_1 - y_2|,$$
so that 
\begin{equation}\label{XX}
   x_1 \asymp x_2 \asymp X := Y_0^{1/2} V^{1/2}.
 \end{equation}  
 Thus we consider now individually the various pieces $\Phi^{\pm}_4(y, (Y_0, Y_0), V)$, where $ 1 - u \asymp V$ and the $\pm$ sign in \eqref{small} is prescribed. 
    If $X \leq T^{1-\varepsilon}$, we can derive a slightly sharper version of \eqref{small}. In this case we can insert the uniform asymptotic expansion \eqref{unifasymp} for the Bessel-$K$-function, and  as before we see by Lemma \ref{int-simple} that the $\rho$-integral is negligible  unless
$$\text{arccosh}\left(\frac{\rho}{x_1}\right) \pm \text{arccosh}\left(\frac{\rho}{x_2}\right) \ll T^{\varepsilon - 1}$$
somewhere in $D_{\rho}$ (Notice that our current assumption $X \leq T^{1-\varepsilon}$ implies that we are away from the branch points of the arccosh-function). Since $\Re \, \text{arccosh}(a+ib) \geq \text{arccosh}(a)$ for $a, b \in \Bbb{R}$ and by the bound \eqref{small}, 
this is impossible if the $\pm$-sign is positive, and in the opposite case we get by the mean value theorem that 
\begin{equation}\label{small-new}
x_1 - x_2 \ll \frac{X}{T^{1-\varepsilon}}.
\end{equation}
Redefining $\varepsilon$, this bound holds by \eqref{small} and \eqref{bessel-x} also in the case  $X \geq T^{1-\varepsilon}$. The same argument shows that in general we can assume that the $\pm$-sign in \eqref{small} is negative, since otherwise $X \ll T^{\varepsilon}$, and using the uniform asymptotic expansion \eqref{unifasymp} we see again that the $\rho$-integral would be negligible. 

Let us temporarily consider values of $y_1, y_2$ where  $\Delta \geq Y_0 T^{2\varepsilon - 1}$. Then by \eqref{usize} we have on the one hand 
\begin{equation}\label{concl1}
  V  \asymp |1 - u| \asymp \Delta/Y_0 + O(T^{\varepsilon - 1})  \asymp \Delta/Y_0,
  \end{equation}
on the other hand we have
\begin{equation}\label{concl2}
\begin{split}
|x_1 - x_2| & \asymp |u y_1^{1/2} - y_2^{1/2}|V^{1/2} = y_2^{1/6}|y_1^{1/3} - y_2^{1/3}| V^{1/2} + O(Y_0^{1/2} V^{1/2} T^{\varepsilon- 1})\\
 & \asymp \frac{\Delta}{Y_0^{1/2}} V^{1/2} + O(Y_0^{1/2} V^{1/2} T^{\varepsilon- 1})  \asymp \frac{\Delta}{Y_0^{1/2}} V^{1/2}.
\end{split}
\end{equation}
Combining \eqref{XX}, \eqref{small-new} and  \eqref{concl2}, we obtain 
$$\Delta \asymp \frac{|x_1 - x_2|  Y_0^{1/2}}{V^{1/2}} \ll \frac{X Y_0^{1/2}}{V^{1/2} T^{1-\varepsilon}} = \frac{Y_0}{T^{1-\varepsilon}},$$
a contradiction. So far we have shown that $\Phi_4^{\pm}(y, (Y_0, Y_0), V)$ is negligible unless the $\pm$ sign is negative and  
$\Delta \leq Y_0 T^{2\varepsilon - 1}$. 
We compute as in \eqref{check} that 
\begin{displaymath}
\begin{split}
     y^j & \frac{\partial^j}{\partial y^j} \Bigl( \exp\Bigl(\pm ( ix_1 \sinh v -  i x_2\sinh (v+\eta))\Bigr)  \left(1 + \frac{z}{y}\right)^{\frac{1}{2}i\tau_2} \Bigr)\\
&   \ll_j \left( \sqrt{1-u^2} |\sinh v| \Bigl(\frac{|\tau_2| |z|}{y} + \frac{|z|}{y^{1/2}} + (|u-1| + |\eta|)y^{1/2}\Bigr)\right) \ll T^{\varepsilon}
   \end{split}
   \end{displaymath}
uniformly in $$x_1, x_2 \ll T, \quad z \ll \frac{Y_0}{T^{1-\varepsilon}}, \quad \sqrt{1-u^2} \sinh v \ll \sqrt{1-u^2} \frac{T}{x_1} \ll \frac{T}{Y_0^{1/2}} \ll T^{\varepsilon}, \quad \eta \ll T^{\varepsilon - 1},$$ 
and conclude that we are in situation \eqref{2}. \\ 

\textbf{The case $j=2$.} Since $u \geq 1$, we have necessarily $Y_2 \gg Y_1$. We insert the uniform asymptotic formula \eqref{unifasymp1}, giving a phase
$$\phi(u, \rho; y_1, y_2) =  \epsilon_1  \tilde{\omega}(x_2, \rho)  + \epsilon_2   \tilde{\omega}(x_1, \rho) + 3 \tau_2 \log u + \frac{ \tau_2}{2} \log \frac{y_1}{y_2},$$
and we compute
\begin{equation}\label{rho1}
  \frac{\partial}{\partial \rho }\phi(u, \rho; y_1, y_2) = \epsilon_1 \text{arcsinh}\left(\frac{\rho}{x_2}\right) + \epsilon_2 \text{arcsinh}\left(\frac{\rho}{x_1}\right).   
\end{equation}
Let us first assume  $Y_2 \geq CY_1$ for some  sufficiently large constant $C$,  so that $u \geq 10$, say, and we are away from the branch-points of the square-roots in $x_1$ and $x_2$ and the cut-off point of the integral.   In this case we argue essentially as in the case $j=1$ with exchanged roles of $y_1$ and $y_2$, so we can be brief.  We have $x_2 \asymp Y_2^{1/2}$ and $x_1 \asymp Y_1^{1/3} Y_2^{1/6}$. The $\rho$-integral 
forces $$T^{2-\varepsilon} \ll x_1 \ll x_2,$$ i.e.\ $Y_1^{1/3} Y_2^{1/6} \gg T^{2-\varepsilon}$ under the present assumption $Y_2 \geq CY_1$, and we verify the lower bound for $t_1, t_2$ required for part b) of the lemma as in the case $j=1$, confirming that we are in situation \eqref{3}.  

 We now turn to the  case where $Y_1 \asymp Y_2 \asymp Y_0$, say. 
 As in the proof of the case $j=4$  we put
 $V \asymp |u - 1| \ll 1$, $\Delta := |y_2 - y_1| \ll Y_0,$ so that \eqref{XX} holds. From now on we consider the pieces $\Phi_2^{\pm}(y, (Y_0, Y_0), V)$ individually, where $V \asymp |1-u|$ and $\epsilon_1 = \pm \epsilon_2$. 
 We think of the $u$-range as smoothly extracted, and it has length $\min(V, T^{\varepsilon-1})$ (the latter by \eqref{usize}). 
 
Let us first assume that $\epsilon_1 = \epsilon_2 = \epsilon$ (say). Then  as before the $\rho$-integral forces $X \gg T^{2-\varepsilon/2}$. Moreover, under this condition we have
$$\frac{\partial}{\partial u }\phi(u, \rho; y_1, y_2)  = \epsilon \frac{u^2 \sqrt{\rho^2 + x_1^2} + \sqrt{\rho^2 + x_2^2}}{u(1 - u^2)}  + \frac{3\tau_2}{u} \asymp \frac{X}{V},
$$
so that by Lemma \ref{int-simple} with $B = X/V$, $r \asymp \min(V, T^{\varepsilon - 1})$, $R = T^{\varepsilon - 1}$   the $u$-integral is negligible.

 From now on we assume  $\epsilon_1 \not= \epsilon_2.$  Let us first assume $X \geq T^{1+\varepsilon/2}$ (which by \eqref{XX} is a condition on the relative size of $Y_0$ and $V$).  Then we argue as in \eqref{indeed} and compute 
$$\frac{\partial}{\partial y_1} \Bigl(\phi(u, \rho; y_1, y_2)  + t_1 \log y_1 + t_2 \log y_2 \Bigr)= \pm \frac{\sqrt{\rho^2 + x_1^2}}{2y_1} + \frac{t_1}{y_1} + \frac{ \tau_2}{2y_1},  $$
so that by Lemma \ref{int-simple} the double Mellin transform is negligible unless $|t_1| \gg x_1 \asymp X \geq T^{1 +\varepsilon/2}$, and similarly for $t_2$.  Hence we are in situation \eqref{4} and there is nothing to prove. 

Therefore we may now suppose that 
\begin{equation}\label{sizeX}
   X\ll T^{1 + \varepsilon/2}. 
\end{equation}   
   In this case we start by considering values of $y_1$, $y_2$ with   $\Delta \geq T^{2\varepsilon-1} Y_0$. Then as in the proof of the previous case $j=4$ we see that \eqref{concl1} and \eqref{concl2} hold. 
 By \eqref{rho1} and the mean value theorem we see that 
$$ \Bigl|\frac{\partial}{\partial \rho }\phi(u, \rho; y_1, y_2)\Bigr| \asymp \frac{T|x_1 - x_2|}{X(X+T)} 
\asymp \frac{\Delta V^{1/2} T}{Y_0^{1/2} X(X+T)}  \asymp \frac{\Delta T}{Y_0 (X+T)}
\gg  T^{\frac{3}{2}\varepsilon - 1},$$
hence by the now familiar argument we see that the $\rho$-integral is negligible. 
Thus we have shown that $\Phi^{\pm}_2(y, (Y_0, Y_0), V)$ is negligible unless the $\pm$ sign is negative and 
\begin{equation}\label{Delta}
\Delta \leq  T^{2\varepsilon-1} Y_0
\end{equation}
holds, and we may assume $X \leq T^{1+\varepsilon/2}$, for otherwise we are in situation \eqref{4}. 

Under these assumptions we argue as in \eqref{partial-y} -- \eqref{osc}, but this time with $$\alpha = 4\pi \sqrt{u^2 - 1} \asymp V^{1/2} \ll  T^{\varepsilon} \min\Bigl(\frac{T}{Y_0^{1/2}}, \frac{1}{T^{1/2}}\Bigr)$$ 
by \eqref{sizeX} and \eqref{Delta} combined with \eqref{usize}, and obtain that 
\begin{equation*}
 y^j  \frac{\partial^j}{\partial y^j} e^{i \phi(u, \rho; y, y+z)} \ll_j   
  T^{\varepsilon},
 \end{equation*}
uniformly in $|z| \ll \Delta$, so that we are  in situation \eqref{2}. 

This completes the proof of the lemma.\hfill $\square$

\section{The spectral large sieve} 

We are now prepared to complete the proof of Theorem \ref{thm2}. With $h$ as defined in \eqref{defh}, it suffices to estimate
$$\sum_{N \leq n, m \leq 2N} a_n \bar{a}_m \int \overline{A_{\pi}(n, 1)} A_{\pi}(m, 1) \frac{h(\nu_{\pi})}{\mathcal{N}(\pi)} d\pi.$$
Summing this over $T^{\varepsilon}$ choices of $\mu_0$ gives the desired upper bound. 
We apply the Kuznetsov formula. The diagonal term becomes
$$\frac{1}{192 \pi^5} \sum_{N \leq n \leq 2N} |a_n|^2  \int h(\mu) \text{spec}(\mu) d\mu \ll T^5 \| a \|_2^2.$$

We  show next that the terms $\Sigma_4$ and $\Sigma_5$ are negligible: indeed, for $\Sigma_4$ the summation condition with $(n_1, n_2, m_1, m_2) = (m, 1, n, 1)$ implies $D_1 = mD_2^2$, so that the argument of $\Phi_{w_4}$ is $\epsilon n/(m D_2^{3}) \ll 1$. As mentioned at the end of Section \ref{kernel}, in this range the function $\Phi_{w_4}$ is negligible. Similarly, the summation condition for $\Sigma_5$ implies $D_1 = n\delta$, $D_2 = n\delta^2$ for some $\delta \in \Bbb{N}$, so that the argument of $\Phi_{w_5}$ is $\epsilon m/(n \delta^{3})  \ll 1$, and again the contribution is negligible. \\

We are left with the estimation of $\Sigma_6$. We split the $D_1, D_2$-sums smoothly  into $O((\log N)^2)$ (by \eqref{unless}) dyadic ranges $D_j \asymp X_j$.   After Mellin inversion we are left with estimating 
$$ \sum_{\epsilon_1, \epsilon_2 = \pm 1} \int_{\Re s = 0}  \left| \widehat{\Phi}\left(s, \Bigl(\frac{N X_2}{X_1^2}, \frac{N X_1}{X_2^2}\Bigr)\right)\right| \sum_{\substack{D_1 \asymp X_1\\  D_2 \asymp X_2} }  \left| \sum_{n, m \asymp N} \frac{a_n \bar{a}_m }{X_1X_2}  S(\epsilon_2 , \epsilon_1 m, n, 1; D_1, D_2)  n^{-s_1} m^{-s_2}\right|  \frac{ds}{(2\pi i)^2},$$
and our aim is to obtain the upper bound $T^{2+\varepsilon} N  \| a \|_2 \| b \|_2 $ for this expression. 
We also split the $s_1, s_2$ contour into dyadic ranges $|t_1| \asymp T_1$, $|t_2| \asymp T_2$. 
We apply Proposition \ref{hybrid} along with Lemmas \ref{phi1} -- \ref{phi3}. 
We distinguish two cases. 

Let us first assume that $X_1 \asymp X_2\asymp X$, say.  Then we apply Lemmas \ref{phi1} -- \ref{phi3} with $Y_1 \asymp Y_2 = N/X =: Y \gg T^{2-\varepsilon}$ by \eqref{unless}. By Lemma \ref{phi1} we can truncate the $s$-contours at $T_1, T_2 \ll Y^{1/2+\varepsilon}$, and by Lemmas \ref{phi2} and \ref{phi3}a) we have $(S_1S_2)^{1/2} \ll T^{2+\varepsilon}$ in all cases. This shows the desired upper bound
$$T^{2+\varepsilon}  (XY^{1/2} + N) \| a \|_2 \| b \|_2 \ll T^{2+\varepsilon} N  \| a \|_2 \| b \|_2 .$$
Let us now assume $X_2 \geq C X_1$ for a sufficiently large constant $C$ (the case $X_1 \geq CX_2$ being   identical up to changing indices). In this case, we have $$Y_1 \asymp \frac{N X_2}{X_1^2} \gg \frac{NX_1}{X_2^2} \asymp Y_2,$$ 
so that 
$$X_1 \asymp \frac{N}{Y_1^{2/3}Y_2^{1/3}}, \quad X_2 \asymp \frac{N}{Y_1^{1/3} Y_2^{2/3}}.$$ 
Moreover, by Lemma \ref{phi3}a) we see that \eqref{s1s2b} holds in any case.  
Proposition \ref{hybrid} and Lemmas \ref{phi1}, \ref{phi2}, \ref{phi3}a) and b) yield $T_1 \asymp Y_1^{1/2}$ and $T_2 \asymp Y_2^{1/3}Y_1^{1/6}\gg T^{2-\varepsilon}$ and 
the upper bound
\begin{displaymath}
\begin{split}
&T^{\varepsilon} \| a \|_2 \| b \|_2  \min\left( T^2, \frac{T^5(T_1+T_2)^{1/2}}{T_1T_2}\right)\Biggl(\frac{N (Y_1^{1/2} Y_2^{1/3}Y_1^{1/6})^{1/2}}{Y_1^{1/3} Y_2^{2/3}} \\
&\quad +  \left(\frac{N}{Y_1^{2/3} Y_2^{1/3}} Y_1^{1/4} +  \frac{N}{Y_1^{1/3} Y_2^{2/3}} (Y_2^{1/3} Y_1^{1/6})^{1/2} \right)  \frac{\sqrt{N}}{(N/Y_1^{2/3}Y_2^{1/3} )^{1/2}}+ N \Biggr)\\
\end{split}
\end{displaymath}
\begin{displaymath}
\begin{split}
& \ll  T^{\varepsilon}  \| a \|_2 \| b \|_2   \min\left( T^2, \frac{T^5(T_1+T_2)^{1/2}}{T_1T_2}\right) N \left(\frac{1}{Y_2^{1/2}} + \frac{1}{Y_1^{1/12} Y_2^{1/6}}  + \frac{Y_1^{1/12}}{Y_2^{1/3}}+ 1  \right) \\
& \ll T^{\varepsilon} \| a \|_2 \| b \|_2 N  \left(T^2  + \frac{T^5 Y_1^{1/4}}{Y_1^{2/3} Y_2^{1/3}}\frac{Y_1^{1/12}}{Y_2^{1/3}}\right)\ll T^{\varepsilon} \| a \|_2 \| b \|_2 N  \left(T^2 + \frac{T^5}{Y_1^{1/3}Y_2^{2/3}}\right).
\end{split}
\end{displaymath}
This completes the proof of Theorem \ref{thm2}. \hfill $\square$\\

The proof of Corollary \ref{cor3} is very similar to that of \cite[Theorem 4]{Bl}, based on the uniform approximate functional equation \cite[Proposition 1]{BH}. The Dirichlet series of $L(s, \pi \times f)$ is given by
$$\sum_{n, m} \frac{A_{\pi}(n, m) \lambda(n)}{(nm^2)^s},$$
where $\lambda(n)$ are the Hecke eigenvalues of $f$. By the Hecke relations \cite[Theorem 6.4.11]{Go} (and M\"obius inversion) we can separate $n$ and $m$ getting
$$\sum_{n, m} \frac{\lambda(n)}{(nm^2)^s} \sum_{d\mid (n, m)} \mu(d) A_{\pi}\left(\frac{n}{d}, 1\right)A_{\pi}\left(1, \frac{m}{d}\right) = \sum_{m, d} \frac{\mu(d) A_{\pi}(1, m)}{(m^2d^3)^s} \sum_{n} \frac{\lambda(nd) A_{\pi}(n, 1)}{n^s}.$$ 
The local factor at infinity is given by (see \cite[Theorem 12.3.6]{Go})
$$L_{\infty}(s, \pi \times f) = \prod_{j=1}^3\prod_{\pm} \Gamma_{\Bbb{R}}(s - \mu_j \pm i\tau), \quad \Gamma_{\Bbb{R}}(s) = \pi^{-s/2} \Gamma(s/2)$$
(and this holds for cuspidal and non-cuspidal $\pi$). Our assumptions imply $\mu_1, \mu_2, \mu_3 \asymp T$ for all $\pi$ under consideration and $|\tau| \leq T^{1-\varepsilon}$, $|\Im s| \leq T^{1-\varepsilon}$, so that all 6 spectral parameters of this $L$-function are $\asymp T$. In particular,   the conductor is $O(T^6)$ and the length of the approximate functional equation is $O(T^{3+\varepsilon})$. More precisely, \cite[Proposition 1]{BH} and Mellin inversion show as in \cite[p.\ 724]{Bl} that
$$L(1/2 + it, \pi \times f) \ll T^{\varepsilon} \int_{-T^{\varepsilon}}^{T^{\varepsilon}} \Bigl| \sum_{nm^2d^3 \ll T^{3+\varepsilon}} \frac{\mu(d) A_{\pi}(1, m)\lambda(nd) A_{\pi}(n, 1) }{(nm^2d^3)^{1/2 +\varepsilon +  it + iu}} \Bigr| du + O(T^{-10}). $$
By the Cauchy-Schwarz inequality we obtain
$$L(1/2 + it, \pi \times f)^2 \ll T^{\varepsilon}  \sum_{m^2d^3 \ll T^{3+\varepsilon}} \frac{|A_{\pi}(1, m)|^2}{m d^{3/2}}  \int_{-T^{\varepsilon}}^{T^{\varepsilon}} \sum_{m^2d^3 \ll T^{3+\varepsilon}} \frac{1}{md^{3/2}}\Bigl| \sum_{n \ll  \frac{T^{3+\varepsilon}}{m^2d^3}} \frac{\lambda(dn) A_{\pi}(n, 1)}{n^{1/2+\varepsilon + it + iu}}\Bigr|^2 du.$$
The first factor is $O(T^{\varepsilon})$ by Rankin-Selberg theory and upper bounds for $L$-values and their residues at $s=1$ (\cite{Li}). Applying Theorem \ref{thm2} for fixed $m$ and $d$, we obtain
 $$\int_{T\Omega} |L(1/2 + it, \pi \times f)|^2 d\pi \ll_{t, f} T^{\varepsilon} \sum_{m^2 d^3 \ll T^{3+\varepsilon}}  \frac{1}{md^{3/2}} T^5 \sum_{n \ll T^{3+\varepsilon}} \frac{|\lambda(dn)|^2}{n} \ll T^{5+\varepsilon},$$
again by Rankin-Selberg.  This completes the proof.  \hfill $\square$

\end{document}